\pdfoutput=1
\PassOptionsToPackage{unicode}{hyperref}
\PassOptionsToPackage{hyphens}{url}
\PassOptionsToPackage{dvipsnames,svgnames,x11names}{xcolor}
\documentclass[
]{article}
\usepackage{xcolor}
\usepackage[margin=1in]{geometry}
\usepackage{amsmath,amssymb}
\usepackage{iftex}
\ifPDFTeX
  \usepackage[T1]{fontenc}
  \usepackage[utf8]{inputenc}
  \usepackage{textcomp} 
\else 
  \usepackage{unicode-math} 
  \defaultfontfeatures{Scale=MatchLowercase}
  \defaultfontfeatures[\rmfamily]{Ligatures=TeX,Scale=1}
\fi
\usepackage{lmodern}
\ifPDFTeX\else
\fi
\IfFileExists{upquote.sty}{\usepackage{upquote}}{}
\IfFileExists{microtype.sty}{
  \usepackage[]{microtype}
  \UseMicrotypeSet[protrusion]{basicmath} 
}{}
\makeatletter
\@ifundefined{KOMAClassName}{
  \IfFileExists{parskip.sty}{%
    \usepackage{parskip}
  }{
    \setlength{\parindent}{0pt}
    \setlength{\parskip}{6pt plus 2pt minus 1pt}}
}{
  \KOMAoptions{parskip=half}}
\makeatother
\makeatletter
\ifx\paragraph\undefined\else
  \let\oldparagraph\paragraph
  \renewcommand{\paragraph}{
    \@ifstar
      \xxxParagraphStar
      \xxxParagraphNoStar
  }
  \newcommand{\xxxParagraphStar}[1]{\oldparagraph*{#1}\mbox{}}
  \newcommand{\xxxParagraphNoStar}[1]{\oldparagraph{#1}\mbox{}}
\fi
\ifx\subparagraph\undefined\else
  \let\oldsubparagraph\subparagraph
  \renewcommand{\subparagraph}{
    \@ifstar
      \xxxSubParagraphStar
      \xxxSubParagraphNoStar
  }
  \newcommand{\xxxSubParagraphStar}[1]{\oldsubparagraph*{#1}\mbox{}}
  \newcommand{\xxxSubParagraphNoStar}[1]{\oldsubparagraph{#1}\mbox{}}
\fi
\makeatother

\usepackage{longtable,booktabs,array}
\usepackage{calc} 
\usepackage{etoolbox}
\makeatletter
\patchcmd\longtable{\par}{\if@noskipsec\mbox{}\fi\par}{}{}
\makeatother
\IfFileExists{footnotehyper.sty}{\usepackage{footnotehyper}}{\usepackage{footnote}}
\makesavenoteenv{longtable}
\usepackage{graphicx}
\makeatletter
\newsavebox\pandoc@box
\newcommand*\pandocbounded[1]{
  \sbox\pandoc@box{#1}%
  \Gscale@div\@tempa{\textheight}{\dimexpr\ht\pandoc@box+\dp\pandoc@box\relax}%
  \Gscale@div\@tempb{\linewidth}{\wd\pandoc@box}%
  \ifdim\@tempb\p@<\@tempa\p@\let\@tempa\@tempb\fi
  \ifdim\@tempa\p@<\p@\scalebox{\@tempa}{\usebox\pandoc@box}%
  \else\usebox{\pandoc@box}%
  \fi%
}
\def\fps@figure{htbp}
\makeatother

\NewDocumentCommand\citeproctext{}{}
\NewDocumentCommand\citeproc{mm}{%
  \begingroup\def\citeproctext{#2}\cite{#1}\endgroup}
\makeatletter
 \let\@cite@ofmt\@firstofone
 \def\@biblabel#1{}
 \def\@cite#1#2{{#1\if@tempswa , #2\fi}}
\makeatother
\newlength{\cslhangindent}
\newlength{\csllabelwidth}
\newenvironment{CSLReferences}[2] 
 {\begin{list}{}{%
  \setlength{\itemindent}{0pt}
  \setlength{\leftmargin}{0pt}
  \setlength{\parsep}{0pt}
  \ifodd #1
   \setlength{\leftmargin}{\cslhangindent}
   \setlength{\itemindent}{-1\cslhangindent}
  \fi
  \setlength{\itemsep}{#2\baselineskip}}}
 {\end{list}}
\usepackage{calc}

\providecommand{\tightlist}{%
  \setlength{\itemsep}{0pt}\setlength{\parskip}{0pt}}

\usepackage{fancyhdr}
\makeatletter
\@ifpackageloaded{caption}{}{\usepackage{caption}}
\AtBeginDocument{%
\ifdefined\contentsname
  \renewcommand*\contentsname{Table of contents}
\else
  \newcommand\contentsname{Table of contents}
\fi
\ifdefined\listfigurename
  \renewcommand*\listfigurename{List of Figures}
\else
  \newcommand\listfigurename{List of Figures}
\fi
\ifdefined\listtablename
  \renewcommand*\listtablename{List of Tables}
\else
  \newcommand\listtablename{List of Tables}
\fi
\ifdefined\figurename
  \renewcommand*\figurename{Figure}
\else
  \newcommand\figurename{Figure}
\fi
\ifdefined\tablename
  \renewcommand*\tablename{Table}
\else
  \newcommand\tablename{Table}
\fi
}
\@ifpackageloaded{float}{}{\usepackage{float}}
\floatstyle{ruled}
\@ifundefined{c@chapter}{\newfloat{codelisting}{h}{lop}}{\newfloat{codelisting}{h}{lop}[chapter]}
\floatname{codelisting}{Listing}

\usepackage{amsthm}
\theoremstyle{plain}
\newtheorem{lemma}{Lemma}[section]
\theoremstyle{plain}
\newtheorem{theorem}{Theorem}[section]
\theoremstyle{plain}
\newtheorem{proposition}{Proposition}[section]
\theoremstyle{plain}
\newtheorem{corollary}{Corollary}[section]
\theoremstyle{remark}
\AtBeginDocument{}

\newtheorem{refremark}{Remark}[section]

\makeatother
\makeatletter
\@ifpackageloaded{caption}{}{\usepackage{caption}}
\@ifpackageloaded{subcaption}{}{\usepackage{subcaption}}
\makeatother
\usepackage{bookmark}
\IfFileExists{xurl.sty}{\usepackage{xurl}}{} 
\hypersetup{
  pdftitle={A Stieltjes representation of order p for \textbackslash operatorname\{tr\}{[}(A+tB)\^{}\{-p\}{]}},
  pdfauthor={Luc Rey-Bellet},
  pdfkeywords={Stieltjes transform, complete
monotonicity, Bessis--Moussa--Villani conjecture, Stahl's
theorem, tracial joint spectral measure, B-splines},
  colorlinks=true,
  linkcolor={blue},
  filecolor={Maroon},
  citecolor={blue},
  urlcolor={blue},
  pdfcreator={LaTeX via pandoc}}

\title{A Stieltjes representation of order \(p\) for
\(\operatorname{tr}[(A+tB)^{-p}]\)}
\usepackage{etoolbox}
\makeatletter
\providecommand{\subtitle}[1]{
  \apptocmd{\@title}{\par {\large #1 \par}}{}{}
}
\makeatother
\subtitle{Heinävaara's measure, B-splines, and Stahl's theorem}
\author{Luc Rey-Bellet}
\date{November 2, 2026}
\begin{document}
\maketitle
\begin{abstract}
For \(A\) positive definite, \(B\) positive semidefinite and \(p\ge1\)
we show that
\(\operatorname{tr}[(A+tB)^{-p}]=\int(1+t\ell)^{-p}\,d\rho_p(\ell)\) for
a positive measure \(\rho_p\) on the spectral interval of
\(L=A^{-1/2}BA^{-1/2}\), and we identify the measure: it has atoms
\(\operatorname{tr}[(P_sA^{-1}P_s)^p]\) at the eigenvalues of \(L\) and
a continuous density, an integral of Heinävaara's tracial joint spectral
measure along rays. For integer \(p\) the positivity of \(\rho_p\)
follows from the positivity of Heinävaara's measure alone, and the
density is a sum of Curry--Schoenberg B-splines; for real \(p\) the
representation is obtained by inverting Heinävaara's transform. For
\(0<p<1\) the function is a Stieltjes function of order \(p\) only if
\([A,B]=0\), but a positive representation of order one holds, with an
explicit density. Along \(A_p=pI+A_0\), \(p\to\infty\), the measures
converge to the measure of Stahl's theorem (the Bessis--Moussa--Villani
conjecture) for \(\operatorname{tr}\,e^{-(A_0+tB)}\).
\end{abstract}



\begin{quote}
\emph{For Jürg Fröhlich, on his eightieth birthday}

Jürg Fröhlich guided my first steps in research during my Diplomarbeit
at ETH Zürich and held me to a standard I have striven to live up to
ever since. While my own work has since wandered some distance from
mathematical physics, it has recently found its way back through a paper
on quantum distributionally robust optimization. This note, a companion
to that work, is offered with admiration and gratitude to the teacher
who set that path in motion.
\end{quote}

\section{Introduction}\label{introduction}

For \(A\) positive definite and \(B\) positive semidefinite Hermitian
matrices on \(\mathbb{C}^n\) and \(p > 0\) set \[
F_p(t) := \operatorname{tr}\bigl[(A + tB)^{-p}\bigr], \qquad t \ge 0 .
\] We ask whether \(F_p\) is a \emph{Stieltjes function of order \(p\)},
that is, whether there is a finite positive Borel measure \(\rho_p\)
with \[
F_p(t) = \int_0^\infty (1 + t\ell)^{-p}\,d\rho_p(\ell), \qquad t\ge0,
\] and, when there is, what \(\rho_p\) looks like.

\textbf{The Bessis--Moussa--Villani problem.} In 1975 Bessis, Moussa and
Villani (\citeproc{ref-bessis1975}{Bessis, Moussa, and Villani 1975}),
seeking monotonically convergent Padé approximants to partition
functions \(\operatorname{tr}\,e^{-(A_0+tB)}\) of quantum statistical
mechanics, conjectured that \(t\mapsto\operatorname{tr}\,e^{-(A_0+tB)}\)
is the Laplace transform of a positive measure for \(A_0\) Hermitian and
\(B\succeq0\). They proved it for \(2\times2\) matrices, where they
found the measure explicitly: two atoms and a positive Bessel density.
Lieb and Seiringer (\citeproc{ref-lieb-seiringer2004}{Lieb and Seiringer
2004}) showed the conjecture equivalent to the complete monotonicity of
\(F_p\) for every \(p>0\). Stahl (\citeproc{ref-stahl2013}{Stahl 2013})
proved it by a Riemann-surface argument, streamlined by Eremenko
(\citeproc{ref-eremenko2015}{Eremenko 2015}) (see also
(\citeproc{ref-clivaz2016}{Clivaz 2016})). Heinävaara
(\citeproc{ref-heinavaara2024tracial}{Heinävaara 2025}) introduced the
\emph{tracial joint spectral measure} \(\mu_{A,B}\), a positive measure
on \(\mathbb R^2\) representing \(\operatorname{tr}\,h(xA+yB)\) for the
functions \(h\) in the range of an explicit integral transform, from
which Stahl's theorem follows and which, in his thesis
(\citeproc{ref-heinavaara-thesis}{Heinävaara 2024b}), yields a
manifestly positive formula for Stahl's density.

\textbf{From Laplace to Stieltjes.} In its Lieb--Seiringer form, Stahl's
theorem says that \(F_p\) is the Laplace transform of a positive measure
on \([0,\infty)\). An order-\(p\) Stieltjes representation is strictly
stronger: the Stieltjes functions of order \(p\) form a subclass of the
completely monotone functions Karp and Prilepkina
(\citeproc{ref-karp2012generalized}{2012}). We will show that, for
\(p\ge 1\), the measure lives on \([l_1,l_n]\), the spectral interval of
\(L:=A^{-1/2}BA^{-1/2}\), it has atoms exactly at the eigenvalues of
\(L\), with explicit weights, and its continuous part, non-zero for
\(p>1\) unless \(A\) and \(B\) commute, has an explicit density,
piecewise polynomial for integer \(p\). Our tool is Heinävaara's
measure: \(t^{-p}\) is not admissible in his theorem, but an explicit
truncated inverse of Heinävaara's transform reduces \(F_p\) to
admissible pieces, and integrating \(\mu_{A,B}\) along the rays
\(b=\ell a\) yields the following. Here \(\lambda_1<\cdots<\lambda_r\)
are the distinct eigenvalues of \(L\), \(P_s\) the orthogonal projection
onto the \(\lambda_s\)-eigenspace, \(\mathcal C_L(Y):=\sum_sP_sYP_s\)
the pinching onto the eigenspaces of \(L\), and \(\mu_c\) the continuous
density of \(\mu_{A,B}\) (Theorem~\ref{thm-tracialjointmeasure}).

\begin{theorem}[Stieltjes-\(p\)
representation]\protect\hypertarget{thm-rep}{}\label{thm-rep}

Let \(A\succ0\) and \(B\succeq0\) be Hermitian on \(\mathbb C^n\) and
\(p>0\). There is a \emph{unique} finite signed Borel measure
\(\rho_p=\rho_{A,B,p}\) on \([0,\infty)\) with
\begin{equation}\phantomsection\label{eq-main-rep}{
\operatorname{tr}\bigl[(A+tB)^{-p}\bigr]=\int_{[0,\infty)}(1+t\ell)^{-p}\,d\rho_p(\ell),
\qquad t\ge0 .
}\end{equation} It is supported in \([l_1,l_n]\) and decomposes as
\(\rho_p=\rho_p^{\mathrm{at}}+\rho_p^{\mathrm{ac}}\) with atoms exactly
at the eigenvalues of \(L\),
\begin{equation}\phantomsection\label{eq-atomic-Omega}{
\rho_p^{\mathrm{at}}=\sum_{s=1}^r\Omega_s\,\delta_{\lambda_s},\qquad
\Omega_s:=\operatorname{tr}\bigl[(P_sA^{-1}P_s)^{p}\bigr]>0,
}\end{equation} and an absolutely continuous part
\begin{equation}\phantomsection\label{eq-cts}{
d\rho_p^{\mathrm{ac}}(\ell)=p(p-1)\Bigl[\int_0^\infty a^{1-p}\mu_c(a,\ell a)\,da\Bigr]d\ell .
}\end{equation} Its total mass and first moment are
\(\operatorname{tr}(A^{-p})\) and \(\operatorname{tr}(A^{-p-1}B)\), and
its continuous mass is
\begin{equation}\phantomsection\label{eq-mass-identity}{
\rho_p^{\mathrm{ac}}\bigl([l_1,l_n]\bigr)=\operatorname{tr}(A^{-p})-\operatorname{tr}\bigl[(\mathcal C_LA^{-1})^{p}\bigr],
}\end{equation} of sign \(\operatorname{sgn}(p-1)\). Consequently \[
\text{(i) }\ \rho_p\ge0\ \text{ for }p\ge1;
\qquad
\text{(ii) for }0<p<1,\ \ \rho_p\ge0\iff[A,B]=0;
\] and in both ranges \(\rho_p^{\mathrm{ac}}=0\) iff \(p=1\) or
\([A,B]=0\). In particular, for \(0<p<1\) the function \(F_p\) is a
Stieltjes function of order \(p\), that is, of the form
(\ref{eq-main-rep}) with a positive measure, if and only if \([A,B]=0\).

\end{theorem}

Two special cases orient the statement. If \([A,B]=0\), simultaneous
diagonalization gives \(\rho_p=\sum_ja_j^{-p}\delta_{b_j/a_j}\), with
\(a_j,b_j\) the eigenvalues of \(A,B\) on a common eigenvector, purely
atomic. At \(p=1\), trace cyclicity gives
\(F_1(t)=\operatorname{tr}[(I+tL)^{-1}A^{-1}]=\sum_iN_i/(1+tl_i)\) with
\(N_i:=\langle e_i|A^{-1}|e_i\rangle\) in an \(L\)-eigenbasis, so
\(\rho_1=\sum_iN_i\delta_{l_i}\) is again atomic. The theorem says these
are the only two ways the continuous part can vanish: away from \(p=1\)
it appears exactly when \(A\) and \(B\) fail to commute, with positive
mass for \(p>1\) and negative mass for \(p<1\). For simple eigenvalues
the atoms are \(N_i^p\), the contribution of the singular rays of
Heinävaara's measure. At a repeated eigenvalue the weight
\(\operatorname{tr}[(P_sA^{-1}P_s)^p]\) differs in general from
\(\sum N_i^p\). In its existence part the theorem follows from
Heinävaara's published theorem by the scheme of his thesis
(\citeproc{ref-heinavaara-thesis}{Heinävaara 2024b}, Remark 2 and proof
of Thm. 24), using a truncation of \(f(x)=p(p-1)x^{-p}\), the preimage
of \(x^{-p}\) under Heinävaara's transform (or alternatively a gamma
integral). Repeated eigenvalues are reached from simple ones by
perturbation and uniqueness. For integer \(p\) the positivity of
\(\rho_p\) (but not its explicit form) is equivalent to the case of a
positive semidefinite direction in Heinävaara's theorem on derivatives
of trace functions (\citeproc{ref-heinavaara2024tracial}{Heinävaara
2025}, Cor. 1.2), for derivatives of order \(p\)
(Section~\ref{sec-bsplines}). For non-integer \(p\) we know no proof
other than the explicit one given here.

The note is organized as follows. Section~\ref{sec-heinavaara} recalls
Heinävaara's theorem and the two arguments he draws from it: a soft one,
which uses only the positivity of \(\mu_{A,B}\) and gives sign
information for the derivatives of trace functions, and an explicit one,
which uses the form of \(\mu_{A,B}\) and gives representing measures.
Section~\ref{sec-bsplines} treats integer \(p\), where everything is
classical or Heinävaara's: \(\rho_p\) is the measure representing the
\(p\)-th derivative of \(s\mapsto\operatorname{tr}f(L+sA^{-1})\), a sum
of Curry--Schoenberg B-splines with weights that are cyclic products of
entries of \(A^{-1}\) (Theorem~\ref{thm-bspline}), and its positivity is
the soft argument. Section~\ref{sec-construction} proves
Theorem~\ref{thm-rep} for all real positive \(p\) by the explicit
argument. For \(n=2\) everything is explicit: the finite-\(p\) density
is a Gauss hypergeometric function (Section~\ref{sec-n2}), and its
\(p\to\infty\) confluence is the Bessel-\(I_1\) density of Bessis,
Moussa and Villani (Section~\ref{sec-n2-stahl}). Along \(A_p=pI+A_0\)
the rescaled measures converge to Stahl's measure of
\(\operatorname{tr}\,e^{-(A_0+tB)}\), whose density is thereby recovered
as a limit of the finite-\(p\) family (Theorem~\ref{thm-stahl-limit},
Section~\ref{sec-stahl}). For \(0<p<1\) the order-\(p\) measure is
signed unless \(A\) and \(B\) commute, but, by operator monotonicity, a
Stieltjes representation of order one on \([0,l_n]\) survives, with an
explicit density (Proposition~\ref{prp-stieltjes1},
Theorem~\ref{thm-rho1-density}, Section~\ref{sec-pless1}). The
\hyperref[sec-appendix-dro]{Appendix} explains the connection with the
companion paper (\citeproc{ref-reybellet2026dro}{Rey-Bellet 2026}).

\subsection{Motivation}\label{sec-motivation}

The problem arose in the companion paper
(\citeproc{ref-reybellet2026dro}{Rey-Bellet 2026}) on distributional
robustness for quantum states. For the sandwiched Rényi divergence of
order \(\alpha\in[\tfrac12,1)\) the dual formula of that paper involves
the trace function \(\operatorname{tr}[(A+tB)^{-p}]\), in the pencil
form of Corollary~\ref{cor-pencil}, with \(p=\alpha/(1-\alpha)\), and
the Stieltjes-\(p\) representation of Theorem~\ref{thm-rep} is what
turns the quantum DRO problem into a classical one. The
\hyperref[sec-appendix-dro]{Appendix} gives the dictionary and gathers
the statements used there. One coincidence is worth recording: under
\(p=\alpha/(1-\alpha)\) the sign change of \(\rho_p\) at \(p=1\) is
exactly \(\alpha=\tfrac12\), the threshold from which the sandwiched
divergence satisfies the data-processing inequality
(\citeproc{ref-frank-lieb2013}{Frank and Lieb 2013}) (and below which it
fails (\citeproc{ref-berta2017}{Berta, Fawzi, and Tomamichel 2017})) and
its Hessian metric is monotone
(\citeproc{ref-takahashi-fujiwara2017}{Takahashi and Fujiwara 2017}):
three thresholds, from three different theories, coincide.

\section{Heinävaara's tracial joint spectral
measure}\label{sec-heinavaara}

\subsection{Heinävaara's theorem}\label{heinuxe4vaaras-theorem}

We recall the result of (\citeproc{ref-heinavaara2024tracial}{Heinävaara
2025}) on which everything rests.

\begin{theorem}[Tracial joint spectral measure
(\citeproc{ref-heinavaara2024tracial}{Heinävaara 2025}, Thms. 1.3 and
1.6)]\protect\hypertarget{thm-tracialjointmeasure}{}\label{thm-tracialjointmeasure}

Let \(A,B\) be Hermitian on \(\mathbb C^n\).

\begin{enumerate}
\def\labelenumi{(\roman{enumi})}
\item
  There is a positive measure \(\mu_{A,B}\) on \(\mathbb R^2\), unique
  away from the origin (\citeproc{ref-heinavaara2024tracial}{Heinävaara
  2025}, Prop. 3.9) and normalized by \(\mu_{A,B}(\{0\})=0\), such that,
  for every measurable \(f:\mathbb R\to\mathbb R\) with \(f(0)=0\) and
  \(\int_{-M}^{M}|f(t)/t|\,dt<\infty\) for all \(M>0\), and every
  \((x,y)\in\mathbb R^2\),
  \begin{equation}\phantomsection\label{eq-heinavaara}{
  \operatorname{tr}\,H(f)(xA+yB)=\int_{\mathbb R^2}f(xa+yb)\,d\mu_{A,B}(a,b),
  \qquad H(f)(z):=\int_0^1\frac{1-t}{t}\,f(zt)\,dt .
  }\end{equation}
\item
  Let \(\mu_{A,B}=\mu_c+\mu_s\) be the Lebesgue decomposition. If the
  pencil \((A,B)\) is non-degenerate (some real combination is
  invertible), \(\mu_c\) has, for almost every \((a,b)\) where \(bA-aB\)
  is invertible, the density
  \begin{equation}\phantomsection\label{eq-rho-c}{
  \mu_c(a,b)=\frac{1}{2\pi}\sum_{i=1}^n\bigl|\operatorname{Im}\lambda_i(C_1(a,b))\bigr|,
  \qquad
  C_1(a,b)=\Bigl(I-\frac{aA+bB}{a^2+b^2}\Bigr)(bA-aB)^{-1},
  }\end{equation} and if moreover \(A\) is invertible and \(A^{-1}B\)
  has \(n\) distinct eigenvalues, then for
  \(\varphi\in C_c^\infty(\mathbb R^2\setminus\{0\})\)
  \begin{equation}\phantomsection\label{eq-mu-s}{
  \mu_s(\varphi)=\sum_{v\in E(A^{-1}B)}\int_0^1\frac{1-t}{t}\,
  \varphi\bigl(t\langle v|A|v\rangle,\,t\langle v|B|v\rangle\bigr)\,dt,
  }\end{equation} \(E(A^{-1}B)\) containing one unit eigenvector per
  eigenvalue.
\end{enumerate}

\end{theorem}

We write \(\mu_c\) both for the continuous part and for its density.

Besides the identity itself, we use three facts about the support and
the size of \(\mu_{A,B}\), all of them in Heinävaara's paper.

\begin{lemma}[Support and finiteness
(\citeproc{ref-heinavaara2024tracial}{Heinävaara 2025}, Lemma 1.5, Prop.
3.8 and proof of Prop.
3.9)]\protect\hypertarget{lem-cone}{}\label{lem-cone}

Let \(A,B\) be Hermitian. Then:

\begin{enumerate}
\def\labelenumi{(\roman{enumi})}
\item
  \emph{(finiteness)} the measure \(\mu_{A,B}\) is in general not
  finite, but \(\widetilde\mu:=(a^2+b^2)\,\mu_{A,B}\) is, with total
  mass
  \(\tfrac16\bigl(\operatorname{tr}A^2+\operatorname{tr}B^2\bigr)\);
\item
  \emph{(half-planes)} if \(xA+yB\succeq0\) for some real \((x,y)\),
  then \(\mu_{A,B}\) is carried by the half-plane \(\{xa+yb\ge0\}\). In
  particular, for \(A\succ0\) and \(B\) Hermitian, with \(l_1\le l_n\)
  the extreme eigenvalues of \(L=A^{-1/2}BA^{-1/2}\),
  \begin{equation}\phantomsection\label{eq-cone}{
  \operatorname{supp}\mu_{A,B}\subset\bigl\{(a,b):\ a\ge0,\ \ l_1a\le b\le l_na\bigr\},
  }\end{equation} so that \(\mu_{A,B}\) is carried by rays of slope
  \(b/a\in[l_1,l_n]\), the spectral interval of \(L\). If moreover
  \(B\succeq0\), then \(b\ge0\) there;
\item
  \emph{(support bound)}
  \(\operatorname{supp}\mu_{A,B}\subset\{(a,b):|a|+|b|\le\|A\|+\|B\|\}\).
\end{enumerate}

\end{lemma}

Items (iii) and the finiteness in (i) are in the proof of
(\citeproc{ref-heinavaara2024tracial}{Heinävaara 2025}, Prop. 3.9). The
mass is (\ref{eq-heinavaara}) with \(f(t)=t^2\), for which
\(H(f)(z)=z^2/6\), at \((x,y)=(1,0)\) and \((0,1)\). Item (ii) is
(\citeproc{ref-heinavaara2024tracial}{Heinävaara 2025}, Lemma 1.5),
stated there for \((x,y)=(1,0)\). The general direction follows from the
change of basis (\citeproc{ref-heinavaara2024tracial}{Heinävaara 2025},
Prop. 3.8), and the cone (\ref{eq-cone}) from the three positive
semidefinite matrices \(A\), \(B-l_1A\) and \(l_nA-B\).

\textbf{The singular part.} When the eigenvalues of \(L\) are
\emph{distinct}, (\ref{eq-mu-s}) puts the singular part on the rays of
slope \(l_i\), with an explicit profile: if \(l_i\) is a simple
eigenvalue of \(L\) with normalized eigenvector \(e_i\), then
\(v_i=A^{-1/2}e_i/\sqrt{N_i}\) (\(N_i=\langle e_i|A^{-1}|e_i\rangle\))
is a normalized eigenvector of \(A^{-1}B\) with
\(\langle v_i|A|v_i\rangle=1/N_i\) and
\(\langle v_i|B|v_i\rangle=l_i/N_i\), and
\begin{equation}\phantomsection\label{eq-singpart}{
\mu_s=\sum_{i=1}^n\int_0^1\frac{1-s}{s}\,\delta_{(s/N_i,\,sl_i/N_i)}(da,db)\,ds
\qquad(\text{distinct eigenvalues}).
}\end{equation} Each ray carries a segment, from the origin to the point
\((1/N_i,\,l_i/N_i)\), with a density that vanishes at the far end and
is not integrable at the origin. At a repeated eigenvalue no such
formula is available, and the proof of Theorem~\ref{thm-rep} does not
need one.

\textbf{The continuous part.} The density (\ref{eq-rho-c}) is explicit
at every multiplicity. Heinävaara's thesis shows that it is supported in
the joint numerical range of the pair,
\(\operatorname{supp}\mu_c\subset W(A,B)=\{(\langle v|A|v\rangle,\langle v|B|v\rangle):\|v\|=1\}\)
(\citeproc{ref-heinavaara-thesis}{Heinävaara 2024b}, Prop. 4). In
particular, \(a\ge\lambda_{\min}(A)\) on \(\operatorname{supp}\mu_c\)
when \(A\succ0\). We shall not need the thesis for this: the consequence
we use comes out of the proof of Theorem~\ref{thm-rep}
(\ref{eq-support-bound}).

\begin{figure}[tb]

\centering{

\includegraphics[width=0.58\linewidth,height=\textheight,keepaspectratio]{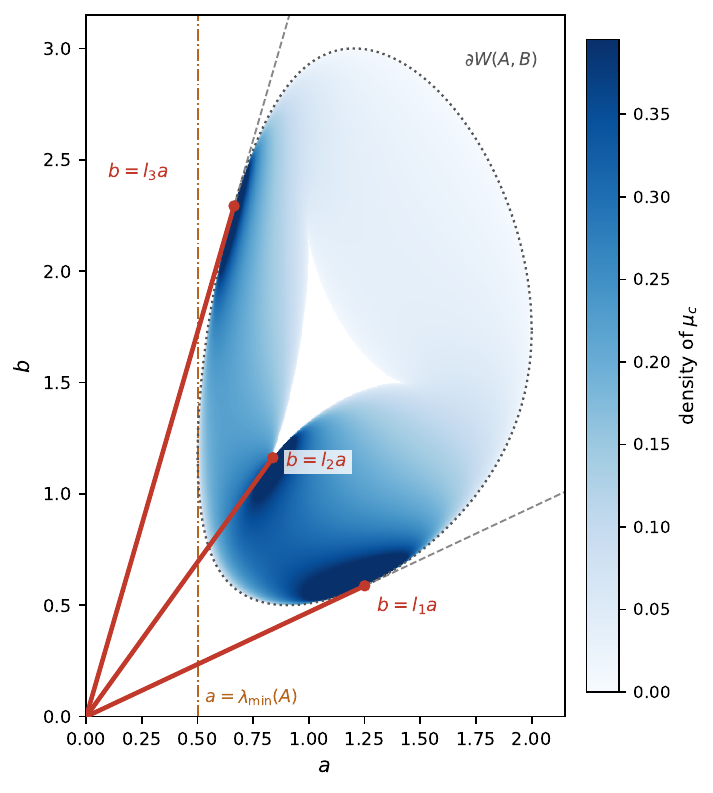}

}

\caption{\label{fig-mu-support}Heinävaara's measure \(\mu_{A,B}\) for
the pair of Figure~\ref{fig-bspline-panels} (\(n=3\), \(A\succ0\),
\(B\succeq0\), \(\operatorname{spec}L\approx\{0.47,1.39,3.45\}\)). The
singular part lives on three segments along the rays \(b=l_ia\), ending
at the points \((1/N_i,l_i/N_i)\). The continuous part (shaded, with the
density (\ref{eq-rho-c})) lies inside the cone (\ref{eq-cone}) between
the extreme rays (dashed), to the right of the line
\(a=\lambda_{\min}(A)\) (\ref{eq-support-bound}), and inside the joint
numerical range \(W(A,B)\) (dotted boundary).}

\end{figure}

Figure~\ref{fig-mu-support} shows the measure for a pair of \(3\times3\)
matrices. The representing measure \(\rho_p\) of Theorem~\ref{thm-rep}
is obtained from it by integrating along the rays through the origin:
the three singular segments give the atoms, and the continuous part
between them gives the density.

\subsection{Two arguments}\label{sec-two-arguments}

Heinävaara draws two kinds of consequence from
Theorem~\ref{thm-tracialjointmeasure}, both of which we will use.

\emph{The soft argument} uses only that \(\mu_{A,B}\) exists, is
positive, and obeys Lemma~\ref{lem-cone}(ii). Applied to the homogeneous
functions \(f(t)=t_+^{k-1}\), for which \(H(f)=f/(k(k-1))\), and
combined with the fact that a function with nonnegative \(k\)-th
derivative is a limit of polynomials of degree less than \(k\) plus
positive combinations of truncated powers \((t-c)_+^{k-1}\)
(\citeproc{ref-bullen1971}{Bullen 1971}, Cor. 8), it gives the
following.

\begin{theorem}[Derivatives of trace functions
(\citeproc{ref-heinavaara2024tracial}{Heinävaara 2025}, Cor.
1.2)]\protect\hypertarget{thm-ktone}{}\label{thm-ktone}

Let \(A,B\) be Hermitian on \(\mathbb C^n\), let \(k\ge1\) be an
integer, and let \(f:\mathbb R\to\mathbb R\) be smooth with
\(f^{(k)}\ge0\). If \(k\) is even, or if \(A\succeq0\), then
\(t\mapsto\operatorname{tr}f(tA+B)\) has nonnegative \(k\)-th
derivative.

\end{theorem}

The cases \(k=1,2\) are classical. The cases \(k=3,4\) were proved
earlier by Heinävaara, by explicit sums of squares, in
(\citeproc{ref-heinavaara2024planes}{Heinävaara 2024a}), where the
general statement was conjectured. The statement has been extended to
infinite dimensions, with the positivity of higher-order spectral shift
functions as an application, by Pradhan and Skripka
(\citeproc{ref-pradhan-skripka2025}{Pradhan and Skripka 2025}). With
\(f(t)=e^{-t}\) and all \(k\) it is Stahl's theorem. With
\(f(t)=t^{-p}\), \(p>0\) (extended, for each \(k\), from a neighbourhood
of the spectra to a smooth function on \(\mathbb R\) whose \(k\)-th
derivative keeps the sign \((-1)^k\), where odd \(k\) uses
\(B\succeq0\)), it gives directly the complete monotonicity of \(F_p\),
statement (iii) of Lieb and Seiringer
(\citeproc{ref-lieb-seiringer2004}{Lieb and Seiringer 2004}, Thm. 1).
The argument says nothing about the representing measures.

\emph{The explicit argument} uses the form of \(\mu_{A,B}\) given in
Theorem~\ref{thm-tracialjointmeasure}(ii), and it produces measures. It
rests on inverting the transform \(H\): if \(h=H(f)\), then
\begin{equation}\phantomsection\label{eq-H-inverse}{
f(x)=2xh'(x)+x^2h''(x)=x\,\bigl(xh(x)\bigr)''
}\end{equation} (\citeproc{ref-heinavaara-thesis}{Heinävaara 2024b},
Remark 2). In his thesis Heinävaara applies (\ref{eq-heinavaara}) at
\((x,y)=(1,0)\), splits \(\mu_{A,B}\) into its singular rays and its
continuous part, and chooses \(f\) by (\ref{eq-H-inverse}) so that
\(H(f)(1/x)=e^{-x}\). The outcome is a manifestly positive formula for
Stahl's measure (\citeproc{ref-heinavaara-thesis}{Heinävaara 2024b, sec.
1.4}, proof of Thm. 24). Section~\ref{sec-construction} applies a
similar scheme to \(x^{-p}\).

For integer \(p\) the two arguments meet: the soft one gives the
positivity of \(\rho_p\) (Section~\ref{sec-bsplines}), the explicit one
its density as an integral of Heinävaara's. For non-integer \(p\) there
is no derivative of order \(p\), and only the explicit argument is
available.

\section{\texorpdfstring{Integer \(p\): the derivative
measure}{Integer p: the derivative measure}}\label{sec-bsplines}

For integer \(p\) the measure \(\rho_p\) of Theorem~\ref{thm-rep} can be
written down, and shown to be positive, without the explicit form of
\(\mu_{A,B}\). Nothing in this section is new except the arrangement.
The identity (\ref{eq-derivative-measure}) below is the Daleckii--Krein
formula for the derivatives of a trace function, in the form used by
Heinävaara (\citeproc{ref-heinavaara2024planes}{Heinävaara 2024a, sec.
2.2}; cf. \citeproc{ref-hiai2010}{Hiai 2010}, Thm. 2.3.1), with the
divided differences written through their Peano kernels
(\citeproc{ref-curry-schoenberg1966}{Curry and Schoenberg 1966, (1.5)}).
In a semifinite von Neumann algebra it is
(\citeproc{ref-skripka2011}{Skripka 2011}, Thm. 3.12). That the measure
it defines is the one representing \(F_p\) is, in moment form, the
duality of Lieb and Seiringer (\citeproc{ref-lieb-seiringer2004}{Lieb
and Seiringer 2004}, Lemma 1) between the pairs \((A,B)\) and
\((L,A^{-1})\) (Remark~\ref{rem-moments}). Its positivity is
Theorem~\ref{thm-ktone}.

Throughout, \(l_1\le\cdots\le l_n\) are the eigenvalues of \(L\) and
\(N_{ij}=\langle e_i|A^{-1}|e_j\rangle\) the entries of \(N:=A^{-1}\) in
an \(L\)-eigenbasis, so that \(N_i=N_{ii}\). For knots
\(\xi_1,\dots,\xi_p\) (repetitions allowed), let
\(M(\cdot\mid\xi_1,\dots,\xi_p)\) be the \emph{Curry--Schoenberg
B-spline}: the Peano kernel of the divided difference
\([\xi_1,\dots,\xi_p]\),
\begin{equation}\phantomsection\label{eq-peano}{
[\xi_1,\dots,\xi_p]\,f=\frac{1}{(p-1)!}\int M(\ell\mid\xi_1,\dots,\xi_p)\,f^{(p-1)}(\ell)\,d\ell ,
}\end{equation} a nonnegative piecewise polynomial of degree \(p-2\)
supported on \([\min_k\xi_k,\max_k\xi_k]\) with \(\int M\,d\ell=1\)
(\citeproc{ref-curry-schoenberg1966}{Curry and Schoenberg 1966,
(1.3)--(1.6)}; \citeproc{ref-deBoor2001}{de Boor 2001};
\citeproc{ref-Karlin1968}{Karlin 1968}). (The kernel appears in this
role already in work of Popoviciu and Tchakaloff in the 1930s, see
(\citeproc{ref-deboor2005divided}{de Boor 2005}).)

\begin{lemma}[B-spline as an order-\(p\) Stieltjes transform
(\citeproc{ref-curry-schoenberg1966}{Curry and Schoenberg 1966}, Lemma
6)]\protect\hypertarget{lem-bspline-pole}{}\label{lem-bspline-pole}

For \(t\ge0\) and knots \(\xi_1,\dots,\xi_p\ge0\), \[
\prod_{k=1}^{p}\frac{1}{1+t\xi_k}=\int_{\min_k\xi_k}^{\max_k\xi_k}(1+t\ell)^{-p}\,M(\ell\mid\xi_1,\dots,\xi_p)\,d\ell .
\] When all knots coincide the left side is \((1+t\xi)^{-p}\) and
\(M(\cdot\mid\xi,\dots,\xi)=\delta_\xi\).

\end{lemma}

This is (\ref{eq-peano}) applied to \(f(\ell)=(c+\ell)^{-1}\) with
\(c=1/t\), for which
\([\xi_1,\dots,\xi_p]f=(-1)^{p-1}\prod_k(c+\xi_k)^{-1}\). The coalescent
case is the limit \(\xi_k\to\xi\), where \(M\to\delta_\xi\).

Expanding the trace in an \(L\)-eigenbasis,
\begin{equation}\phantomsection\label{eq-Bsplines}{
F_p(t)=\operatorname{tr}\bigl[\bigl((I+tL)^{-1}A^{-1}\bigr)^p\bigr]
=\sum_{i_1,\dots,i_p=1}^n\frac{N_{i_1i_2}N_{i_2i_3}\cdots N_{i_pi_1}}{(1+tl_{i_1})\cdots(1+tl_{i_p})},
}\end{equation} and inserting Lemma~\ref{lem-bspline-pole} into each
term gives the representation.

\begin{theorem}[B-spline form of \(\rho_p\), integer
\(p\)]\protect\hypertarget{thm-bspline}{}\label{thm-bspline}

Let \(p\ge1\) be an integer, \(A\succ0\), \(B\succeq0\), and define the
signed measure
\begin{equation}\phantomsection\label{eq-bspline-measure}{
\rho_p:=\sum_{i_1,\dots,i_p=1}^n N_{i_1i_2}N_{i_2i_3}\cdots N_{i_pi_1}\,
M(\cdot\mid l_{i_1},\dots,l_{i_p}) .
}\end{equation}

\begin{enumerate}
\def\labelenumi{(\roman{enumi})}
\item
  \(F_p(t)=\int(1+t\ell)^{-p}\,d\rho_p(\ell)\) for \(t\ge0\), so
  \(\rho_p\) is the measure of Theorem~\ref{thm-rep}.
\item
  \(\rho_p\) represents the \(p\)-th derivative of trace functions in
  the direction \(A^{-1}\): for every \(f\in C^p\) near \([l_1,l_n]\),
  \begin{equation}\phantomsection\label{eq-derivative-measure}{
  \frac{d^p}{ds^p}\operatorname{tr}f\bigl(L+sA^{-1}\bigr)\Big|_{s=0}=\int f^{(p)}(\ell)\,d\rho_p(\ell).
  }\end{equation}
\item
  \(\rho_p\ge0\).
\item
  The atoms of \(\rho_p\) are \(\Omega_s\delta_{\lambda_s}\), with
  \(\Omega_s=\operatorname{tr}[(P_sA^{-1}P_s)^p]\) as in
  Theorem~\ref{thm-rep}, and its continuous part has the density
  \begin{equation}\phantomsection\label{eq-bspline-density}{
  \frac{d\rho_p^{\mathrm{ac}}}{d\ell}(\ell)=\!\!\sum_{\substack{i_1,\dots,i_p\\ l_{i_1},\dots,l_{i_p}\text{ not all equal}}}\!\!
  N_{i_1i_2}\cdots N_{i_pi_1}\,M(\ell\mid l_{i_1},\dots,l_{i_p}),
  }\end{equation} a real piecewise polynomial of degree at most \(p-2\)
  with breakpoints at the eigenvalues of \(L\).
\end{enumerate}

\end{theorem}

\begin{proof}
\leavevmode

\begin{enumerate}
\def\labelenumi{(\roman{enumi})}
\item
  is (\ref{eq-Bsplines}) and Lemma~\ref{lem-bspline-pole}. The
  identification with the measure of Theorem~\ref{thm-rep} is the
  uniqueness statement there, which is proved from the moments alone.
\item
  Both sides are continuous linear functionals of \(f\) in the \(C^p\)
  norm on \([l_1,l_n]\). For the right side this is clear. For the left
  side, the eigenvalues \(\lambda_i(s)\) of \(L+sA^{-1}\) are analytic
  in \(s\) (Rellich), so
  \(\operatorname{tr}f(L+sA^{-1})=\sum_if(\lambda_i(s))\), and its
  \(p\)-th derivative at \(s=0\) is a combination of the numbers
  \(f^{(j)}(l_i)\), \(j\le p\). It therefore suffices to take \(f\)
  analytic near \([l_1,l_n]\), and then, by Cauchy's formula,
  \(f(x)=(z-x)^{-1}\) with \(z\notin\mathbb R\). For this \(f\) the
  Neumann series gives \[
  \frac{d^p}{ds^p}\operatorname{tr}\bigl(zI-L-sA^{-1}\bigr)^{-1}\Big|_{s=0}
  =p!\,\operatorname{tr}\bigl[\bigl((zI-L)^{-1}A^{-1}\bigr)^p(zI-L)^{-1}\bigr]
  =p!\sum_{i_1,\dots,i_p}\frac{N_{i_1i_2}\cdots N_{i_pi_1}}{(z-l_{i_1})\prod_k(z-l_{i_k})} .
  \] On the other hand, Lemma~\ref{lem-bspline-pole} with \(t=-1/z\),
  extended by analyticity in \(t\) to every \(z\notin[l_1,l_n]\), gives
  \(\prod_k(z-\xi_k)^{-1}=\int(z-\ell)^{-p}M(\ell\mid\xi_1,\dots,\xi_p)\,d\ell\),
  and hence
  \(\sum N_{i_1i_2}\cdots N_{i_pi_1}\prod_k(z-l_{i_k})^{-1}=\int(z-\ell)^{-p}\,d\rho_p(\ell)\).
  Differentiating in \(z\), and using the cyclic invariance of the
  weights to collect the \(p\) terms produced on the left, gives
  \(\sum N_{i_1i_2}\cdots N_{i_pi_1}(z-l_{i_1})^{-1}\prod_k(z-l_{i_k})^{-1}=\int(z-\ell)^{-p-1}\,d\rho_p(\ell)\),
  and \(f^{(p)}(\ell)=p!\,(z-\ell)^{-p-1}\).
\item
  Let \(g\ge0\) be smooth with compact support and let \(f\) be a
  \(p\)-fold primitive of \(g\). By Theorem~\ref{thm-ktone} applied to
  the pair \((A^{-1},L)\), with \(k=p\) and \(A^{-1}\succ0\), the
  function \(s\mapsto\operatorname{tr}f(L+sA^{-1})\) has nonnegative
  \(p\)-th derivative, in particular at \(s=0\). By
  (\ref{eq-derivative-measure}), \(\int g\,d\rho_p\ge0\). Such \(g\)
  approximate every nonnegative continuous function uniformly on
  \([l_1,l_n]\), and \(\rho_p\) is a finite signed measure. Hence
  \(\rho_p\ge0\).
\item
  The words with all knots equal to \(\lambda_s\) are those with all
  indices in the index set of the \(\lambda_s\)-eigenspace of \(L\).
  Each contributes its weight times \(\delta_{\lambda_s}\), and the sum
  of these weights is \(\operatorname{tr}[(P_sA^{-1}P_s)^p]\). Every
  other word contributes a B-spline of degree \(p-2\) with knots at
  eigenvalues of \(L\). The sum is real: reversing a cyclic word,
  \((i_1,\dots,i_p)\mapsto(i_1,i_p,\dots,i_2)\), conjugates its weight
  (as \(N_{ji}=\overline{N_{ij}}\)) and fixes its knots.
\end{enumerate}

\end{proof}

The density (\ref{eq-bspline-density}) is constant at \(p=2\), linear at
\(p=3\) and quadratic at \(p=4\) between consecutive eigenvalues. In
general it jumps at the eigenvalues: the words \((i,\dots,i,j)\) carry a
knot of multiplicity \(p-1\) at \(l_i\), where a B-spline of degree
\(p-2\) is discontinuous.

The structure of (\ref{eq-derivative-measure}), words in the
perturbation integrated against splines, also underlies the higher-order
spectral shift functions of perturbation theory
(\citeproc{ref-dykema-skripka2009}{Dykema and Skripka 2009}, Cor. 9.4;
\citeproc{ref-skripka2011}{Skripka 2011}), which are averages of
derivative measures along the segment from \(L\) to \(L+A^{-1}\), and
whose positivity (\citeproc{ref-pradhan-skripka2025}{Pradhan and Skripka
2025}) is deduced from Theorem~\ref{thm-ktone} as in (iii).

\begin{refremark}[Moments and the Lieb--Seiringer duality]
For \(f(x)=x^{p+r}\), (\ref{eq-derivative-measure}) reads
\(\frac{(p+r)!}{r!}\int\ell^r\,d\rho_p=\frac{d^p}{ds^p}\operatorname{tr}(L+sA^{-1})^{p+r}\big|_{s=0}\).
The coefficient of \(s^p\) in \(\operatorname{tr}(L+sA^{-1})^{p+r}\) is
the coefficient of \(\lambda^r\) in
\(\operatorname{tr}(A^{-1}+\lambda L)^{p+r}\), so this is
\begin{equation}\phantomsection\label{eq-ls-duality}{
(p+r)\,\frac{d^r}{dt^r}\operatorname{tr}(A+tB)^{-p}\Big|_{t=0}
=p\,(-1)^r\,\frac{d^r}{d\lambda^r}\operatorname{tr}\bigl(A^{-1}+\lambda L\bigr)^{p+r}\Big|_{\lambda=0},
}\end{equation} the identity by which Lieb and Seiringer pass between
negative powers of one pair and positive powers of the other
(\citeproc{ref-lieb-seiringer2004}{Lieb and Seiringer 2004}, Lemma 1).
They later extended it to all complex \(p\)
(\citeproc{ref-lieb-seiringer2012}{Lieb and Seiringer 2012}, Lemma 1).
For real \(p\) it gives the moments of the measure of
Theorem~\ref{thm-rep}: \(\int\ell^r\,d\rho_p=D_r/(p+1)_r\), with \(D_r\)
the derivative on the right of (\ref{eq-ls-duality}) and
\((p+1)_r=(p+1)\cdots(p+r)\).

\label{rem-moments}

\end{refremark}

\begin{figure}

\centering{

\includegraphics[width=1\linewidth,height=\textheight,keepaspectratio]{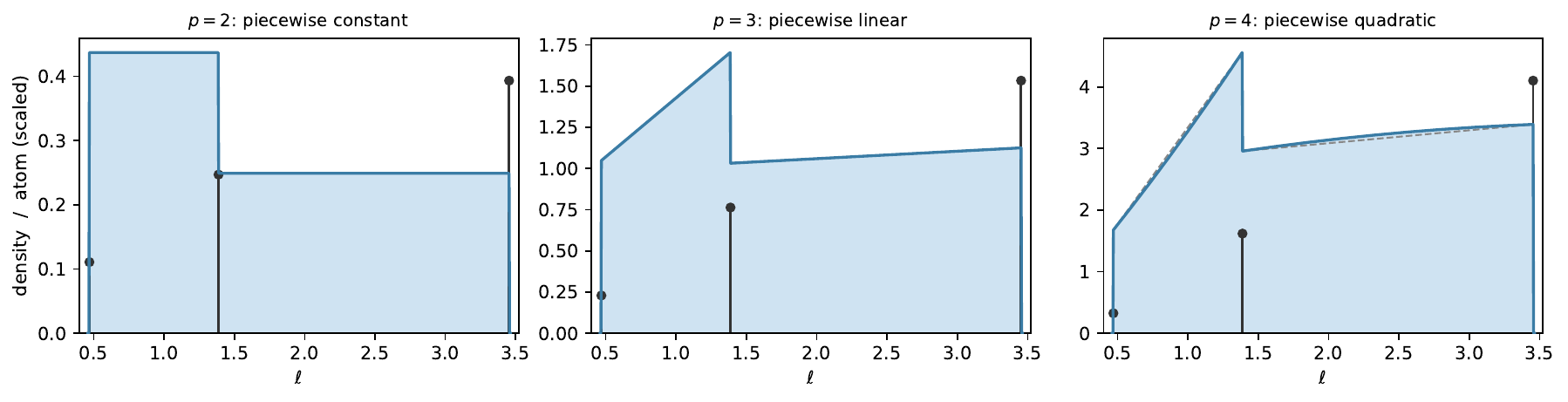}

}

\caption{\label{fig-bspline-panels}The continuous density of \(\rho_p\)
from the B-spline representation (\ref{eq-bspline-measure}), for
\(B=\operatorname{diag}(0.5,1.5,3)\) and \(A\) with eigenvalues
\((0.5,1,2)\) in a rotated basis
(\(\operatorname{spec}L\approx\{0.47,1.39,3.45\}\)): piecewise constant
at \(p=2\), piecewise linear at \(p=3\), piecewise quadratic at \(p=4\),
with breakpoints at \(\operatorname{spec}(L)\). Dashed lines are the
chords across each interval between eigenvalues --- at \(p=3\) the
density lies on them, at \(p=4\) it bows away from them. Stems mark the
atoms, rescaled to the height of the density.}

\end{figure}

Figure~\ref{fig-bspline-panels} shows the density for an \(n=3\) pair at
\(p=2,3,4\).

\textbf{Elementary proofs.} The positivity (iii) rests on
Theorem~\ref{thm-ktone}, and hence on the tracial joint spectral
measure. For small \(p\) it can be seen directly. At \(p=2\) the cyclic
products are \(|N_{ij}|^2\ge0\) and the B-splines are normalized
indicators of \([l_i,l_j]\), so (\ref{eq-bspline-density}) is a
nonnegative combination of boxes. For \(p\ge3\) positivity is not
termwise: the fully off-diagonal products \(N_{ij}N_{jk}N_{ki}\) are
complex, paired reversed words contribute their real part, which may be
negative, and positivity rests on cancellation. For \(p=3\) and \(p=4\)
Heinävaara found the cancellation, as explicit sums of squares valid for
every \(n\) (\citeproc{ref-heinavaara2024planes}{Heinävaara 2024a}). At
\(p=3\) his computation reads as follows in the present notation.

\begin{proposition}[Positivity at \(p=3\), after
Heinävaara]\protect\hypertarget{prp-p3}{}\label{prp-p3}

Let \(p=3\) and let \(\ell\in(l_1,l_n)\) not be an eigenvalue of \(L\).
Write \(N^{>}\) and \(N^{<}\) for the principal submatrices of
\((N_{ij})\) on the indices \(j\) with \(l_j>\ell\) and with
\(l_j<\ell\). Then \begin{equation}\phantomsection\label{eq-p3-sos}{
\begin{aligned}
\frac{d\rho_3^{\mathrm{ac}}}{d\ell}(\ell)
&=6\sum_{i:\,l_i<\ell}(\ell-l_i)\,\langle v_i|N^{>}|v_i\rangle
+6\sum_{i:\,l_i>\ell}(l_i-\ell)\,\langle u_i|N^{<}|u_i\rangle,\\
v_i&:=\Bigl(\frac{N_{ji}}{l_j-l_i}\Bigr)_{l_j>\ell},\qquad
u_i:=\Bigl(\frac{N_{ji}}{l_i-l_j}\Bigr)_{l_j<\ell}.
\end{aligned}
}\end{equation} Every term is nonnegative, because the principal
submatrices of \(A^{-1}\) are positive semidefinite.

\end{proposition}

\begin{proof}
By (\citeproc{ref-curry-schoenberg1966}{Curry and Schoenberg 1966,
(1.3)--(1.4)}),
\(M(\ell\mid\xi_1,\xi_2,\xi_3)=2\,[\xi_1,\xi_2,\xi_3](\cdot-\ell)_+\).
This vanishes when the three knots lie on the same side of \(\ell\). It
equals \(2(\ell-a)/((b_1-a)(b_2-a))\) for knots \(a<\ell<b_1,b_2\), and
\(2(b-\ell)/((b-a_1)(b-a_2))\) for knots \(a_1,a_2<\ell<b\). In
(\ref{eq-bspline-density}) the words with exactly one index \(i\) below
\(\ell\) therefore contribute, the three positions of that index giving
equal sums by cyclicity, \[
3\sum_{l_i<\ell}\ \sum_{l_j,\,l_k>\ell}N_{ij}N_{jk}N_{ki}\,\frac{2(\ell-l_i)}{(l_j-l_i)(l_k-l_i)},
\] which is the first sum in (\ref{eq-p3-sos}). The words with exactly
one index above \(\ell\) give the second.
\end{proof}

At \(p=4\) the corresponding statement is Heinävaara's theorem on fourth
derivatives (\citeproc{ref-heinavaara2024planes}{Heinävaara 2024a}),
which needs no positivity of the direction \(A^{-1}\). For integer
\(p\ge5\) we know of no elementary proof. See
(\citeproc{ref-heinavaara2024planes}{Heinävaara 2024a, sec. 3}) for
\(2\times2\) matrices and partial results.

\section{\texorpdfstring{Real \(p\): the Stieltjes-\(p\)
representation}{Real p: the Stieltjes-p representation}}\label{sec-construction}

We prove Theorem~\ref{thm-rep} for every real \(p>0\), by the explicit
argument of Section~\ref{sec-two-arguments}. The function \(t^{-p}\) is
not admissible in Theorem~\ref{thm-tracialjointmeasure}. We invert
Heinävaara's transform directly, as he does for Stahl's measure
(Section~\ref{sec-two-arguments}). By (\ref{eq-H-inverse}),
\(h(x)=x^{-p}\) corresponds to \(f(x)=p(p-1)x^{-p}\). The singularity at
the origin is handled by a cutoff, whose two boundary terms cancel
exactly against the singular rays of \(\mu_{A,B}\). This is the
mechanism behind both the coefficient \(p(p-1)\) and the exceptional
case \(p=1\). On the continuous part the result is, formally,
Heinävaara's identity
\(\operatorname{tr}|xA+yB|^{q}=q(q+1)\int|ax+by|^{q}\,d\mu_{A,B}\)
(\citeproc{ref-heinavaara2024tracial}{Heinävaara 2025}, after Thm. 1.3)
at \(q=-p\), where his theorem no longer applies. The work is in the
singular rays. A different proof, based on a gamma-integral
representation of \(x^{-p}\), is sketched in
Remark~\ref{rem-gamma-route}. Three lemmas are needed. The first
identifies the atoms at repeated eigenvalues, the second gives the sign
of the continuous mass, and the third passes from simple to repeated
eigenvalues.

\subsection{Three lemmas}\label{sec-singular}

Let \(\lambda_1<\dots<\lambda_r\) be the distinct eigenvalues of \(L\),
with eigenspaces \(V_s\) of dimension \(k_s\) and orthogonal projections
\(P_s\). Set \begin{equation}\phantomsection\label{eq-Ntilde}{
\widetilde N_s:=P_sA^{-1}P_s\big|_{V_s}\succ0\ \text{ on }V_s,\qquad
\nu_{s,1},\dots,\nu_{s,k_s}:=\operatorname{spec}\widetilde N_s,
}\end{equation} and \begin{equation}\phantomsection\label{eq-Phi-def}{
\Phi_p(z):=\operatorname{tr}\bigl[(zA-B)^{-p}\bigr],\qquad \operatorname{Im}z>0 .
}\end{equation} For \(\operatorname{Im}z>0\) and \(w\neq0\),
\(\operatorname{Im}\langle w|zA-B|w\rangle=(\operatorname{Im}z)\langle w|A|w\rangle>0\),
so the spectrum of \(zA-B\) lies in the open upper half-plane. The
principal power is well defined and \(\Phi_p\) is holomorphic on
\(\{\operatorname{Im}z>0\}\).

\begin{lemma}[Matrix
pole]\protect\hypertarget{lem-Apole}{}\label{lem-Apole}

~

\begin{enumerate}
\def\labelenumi{(\roman{enumi})}
\tightlist
\item
  For each \(s\) and every \(p>0\), \[
  \Phi_p(\lambda_s+i\varepsilon)=(i\varepsilon)^{-p}\operatorname{tr}\bigl(\widetilde N_s^{\,p}\bigr)+o(\varepsilon^{-p})
  \qquad(\varepsilon\to0^+).
  \]
\item
  Consequently, if for some constants \(\omega_1,\dots,\omega_r\) and
  some finite signed measure \(\eta\) on \([\lambda_1,\lambda_r]\)
  without atoms \begin{equation}\phantomsection\label{eq-stieltjes-p}{
  \Phi_p(z)=\sum_{s=1}^r\omega_s\,(z-\lambda_s)^{-p}+\int(z-\ell)^{-p}\,d\eta(\ell),
  \qquad\operatorname{Im}z>0,
  }\end{equation} then
  \(\omega_s=\operatorname{tr}(\widetilde N_s^{\,p})\) for every \(s\).
\end{enumerate}

\end{lemma}

\begin{proof}
\leavevmode

\begin{enumerate}
\def\labelenumi{(\roman{enumi})}
\tightlist
\item
  Write \(zA-B=R_s+i\varepsilon A\) with \(z=\lambda_s+i\varepsilon\)
  and \(R_s:=\lambda_sA-B=A^{1/2}(\lambda_sI-L)A^{1/2}\) Hermitian. Then
  \(R_sv=0\iff A^{1/2}v\in
  V_s\), so \(\ker R_s=W_s:=A^{-1/2}V_s\) has dimension \(k_s\) and
  \(R_s\) is invertible on \(W_s^\perp\). As \(0\) is a semisimple
  eigenvalue of the Hermitian \(R_s\), the \(k_s\) eigenvalues of
  \(R_s+i\varepsilon A\) bifurcating from \(0\) are
  (\citeproc{ref-kato1995}{Kato 1995}, Thm. II.5.11) \[
  \mu_j(\varepsilon)=i\varepsilon\,\theta_j+o(\varepsilon),\qquad
  \theta_j\in\operatorname{spec}\bigl(\Pi_sA|_{W_s}\bigr),
  \] with \(\Pi_s\) the orthogonal projection onto \(W_s\), the
  remaining \(n-k_s\) staying bounded away from \(0\). For an
  orthonormal basis \(\{u_a\}\) of \(V_s\) put \(w_a:=A^{-1/2}u_a\).
  Then \(\langle w_a|w_b\rangle=(\widetilde N_s)_{ab}\) and
  \(\langle w_a|A|w_b\rangle=\delta_{ab}\), so
  \(\Pi_sA|_{W_s}=\widetilde N_s^{-1}\) in the basis \(\{w_a\}\), giving
  \(\theta_j=1/\nu_{s,j}\) and
  \(\mu_j(\varepsilon)=(i\varepsilon/\nu_{s,j})(1+o(1))\) in the open
  upper half-plane. Hence \[
  \sum_{j=1}^{k_s}\mu_j(\varepsilon)^{-p}
  =(i\varepsilon)^{-p}\sum_j\nu_{s,j}^{\,p}\,(1+o(1))
  =(i\varepsilon)^{-p}\operatorname{tr}\bigl(\widetilde N_s^{\,p}\bigr)+o(\varepsilon^{-p}),
  \] and the bounded eigenvalues contribute \(O(1)\).
\item
  Fix \(s\) and let \(z=\lambda_s+i\varepsilon\). The term of index
  \(s\) in (\ref{eq-stieltjes-p}) is \(\omega_s(i\varepsilon)^{-p}\).
  The terms of index \(s'\neq s\) are
  \(\omega_{s'}(\lambda_s-\lambda_{s'}+i\varepsilon)^{-p}\), which
  remain bounded as \(\varepsilon\to0^+\). For the integral,
  \(|w^{-p}|=|w|^{-p}\) on the principal branch, so \[
  \varepsilon^{p}\,\bigl|(z-\ell)^{-p}\bigr|
  =\varepsilon^{p}\bigl((\ell-\lambda_s)^2+\varepsilon^2\bigr)^{-p/2}=:g_\varepsilon(\ell).
  \] Here \(0\le g_\varepsilon\le1\) and
  \(g_\varepsilon\to\mathbf 1_{\{\lambda_s\}}\) pointwise as
  \(\varepsilon\to0^+\). Since \(\eta\) is finite, dominated convergence
  gives \[
  \varepsilon^{p}\,\Bigl|\int(z-\ell)^{-p}\,d\eta(\ell)\Bigr|\le\int g_\varepsilon\,d|\eta|
  \longrightarrow|\eta|(\{\lambda_s\})=0,
  \] because \(\eta\) has no atoms. Hence
  \(\Phi_p(\lambda_s+i\varepsilon)=\omega_s(i\varepsilon)^{-p}+o(\varepsilon^{-p})\).
  Subtracting (i),
  \(\bigl(\omega_s-\operatorname{tr}(\widetilde N_s^{\,p})\bigr)(i\varepsilon)^{-p}=o(\varepsilon^{-p})\),
  and since \(|(i\varepsilon)^{-p}|=\varepsilon^{-p}\) this forces
  \(\omega_s=\operatorname{tr}(\widetilde N_s^{\,p})\).
\end{enumerate}

\end{proof}

Let \(N=A^{-1}\) and let \(\mathcal C_L(Y):=\sum_sP_sYP_s\) be the
pinching onto the eigenspaces of \(L\), so that
\(\mathcal C_L(N)=\bigoplus_s\widetilde N_s\).

\begin{lemma}[Pinching
inequality]\protect\hypertarget{lem-pinching}{}\label{lem-pinching}

With \(N=A^{-1}\) and \(\mathcal C_L\) the pinching onto the eigenspaces
of \(L\), for every \(p>0\),
\begin{equation}\phantomsection\label{eq-pinch-ineq}{
\operatorname{tr}\bigl[(\mathcal C_LN)^{p}\bigr]-\operatorname{tr}(N^{p})
\ \begin{cases}\le0,& p>1,\\[0.5mm]=0,& p=1,\\[0.5mm]\ge0,& 0<p<1,\end{cases}
\qquad\text{with equality (for }p\neq1\text{) iff }[N,L]=0,\text{ i.e. }[A,B]=0 .
}\end{equation}

\end{lemma}

\begin{proof}
The pinching is an average of unitary conjugations: with
\(U_\varepsilon=\sum_s\varepsilon_sP_s\)
(\(\varepsilon\in\{\pm1\}^{r}\), self-inverse unitary) and
\(\mathbb E_\varepsilon[\varepsilon_s\varepsilon_{s'}]=\delta_{ss'}\),
\begin{equation}\phantomsection\label{eq-pinch-average}{
\mathcal C_L(N)=2^{-r}\sum_\varepsilon U_\varepsilon NU_\varepsilon^{*}
}\end{equation} (off-diagonal blocks cancel). For pinchings and the
majorization \(\lambda(\mathcal C(Y))\prec\lambda(Y)\) see
(\citeproc{ref-bhatia1997}{Bhatia 1997}, Problem II.5.5). On
\(\{Y\succ0\}\), \(\Psi(Y):=\operatorname{tr}(Y^{p})\) is strictly
convex for \(p>1\) and strictly concave for \(0<p<1\)
(\citeproc{ref-bhatia1997}{Bhatia 1997};
\citeproc{ref-carlen2010}{Carlen 2010}). Since
\(\Psi(U_\varepsilon NU_\varepsilon^{*})=\Psi(N)\), Jensen applied to
(\ref{eq-pinch-average}) gives (\ref{eq-pinch-ineq}), with equality (for
\(p\neq1\)) iff all \(U_\varepsilon NU_\varepsilon^{*}\) coincide,
i.e.~iff \(N\) commutes with every \(P_s\), i.e.~\([A,B]=0\). For
\(p=1\) both sides equal \(\operatorname{tr}(N)\).
\end{proof}

The theorem is proved first for simple eigenvalues, and then extended by
approximating \(B\) by matrices \(B_m\) for which
\(L_m=A^{-1/2}B_mA^{-1/2}\) has simple spectrum. The next lemma controls
the continuous densities in this limit: away from the spectrum of \(L\)
they are uniformly bounded and converge. It says nothing at the
eigenvalues of \(L\) themselves. When several eigenvalues of \(L_m\)
merge, the density between them may concentrate, and part of the
continuous mass is then absorbed into the atom of the limit. This is why
the atoms at repeated eigenvalues are identified by
Lemma~\ref{lem-Apole} and not as limits of atoms.

In the lemma the right side of (\ref{eq-rho-c}) is taken as the
definition of the function \(\mu_c\) at every point where \(bA-aB\) is
invertible.

\begin{lemma}[Convergence of the radial
densities]\protect\hypertarget{lem-density-limit}{}\label{lem-density-limit}

Let \(A\succ0\) and \(p>0\), let \(B_m\succeq0\) converge to \(B\), and
let \(\mu_c^{(m)}\) and \(\mu_c\) be the densities (\ref{eq-rho-c}) of
\((A,B_m)\) and \((A,B)\). Put \[
g_m(\ell):=p(p-1)\int_0^\infty a^{1-p}\mu_c^{(m)}(a,\ell a)\,da,\qquad
g(\ell):=p(p-1)\int_0^\infty a^{1-p}\mu_c(a,\ell a)\,da .
\] Assume that, for every \(m\), \(\mu_c^{(m)}(a,b)=0\) for almost every
\((a,b)\) with \(0<a<\lambda_{\min}(A)\). Then \(\mu_c(a,b)=0\) for
almost every \((a,b)\) with \(0<a<\lambda_{\min}(A)\), and for every
compact \(K\subset\mathbb R\setminus\operatorname{spec}L\) there is
\(m_0\) such that \[
\sup_{m\ge m_0}\,\sup_{\ell\in K}|g_m(\ell)|<\infty
\qquad\text{and}\qquad
g_m(\ell)\to g(\ell)\quad(\ell\in K).
\]

\end{lemma}

\begin{proof}
Let \(L_m:=A^{-1/2}B_mA^{-1/2}\), so that \(L_m\to L\), and choose
\(m_0\) such that the distance from \(K\) to \(\operatorname{spec}L_m\)
is bounded below for \(m\ge m_0\). For \(\ell\in K\), \(a>0\) and
\(b=\ell a\), \[
(bA-aB_m)^{-1}=a^{-1}A^{-1/2}(\ell I-L_m)^{-1}A^{-1/2},
\] and the same holds for \(\ell\) in an open set \(U\supset K\) whose
distance to \(\operatorname{spec}L_m\) is bounded below for
\(m\ge m_0\). So, for \(m\ge m_0\), the function
\((a,\ell)\mapsto\mu_c^{(m)}(a,\ell a)\) is continuous on the open set
\((0,\infty)\times U\), and it converges pointwise there to
\(\mu_c(a,\ell a)\) as \(m\to\infty\), because the eigenvalues of a
matrix depend continuously on it. The map \((a,\ell)\mapsto(a,\ell a)\)
is a diffeomorphism of \(\{a>0\}\), so the hypothesis says that this
continuous function vanishes almost everywhere on the open set
\((0,\lambda_{\min}(A))\times U\). Hence it vanishes everywhere on it.
In the same way Lemma~\ref{lem-cone}(iii) gives
\(\mu_c^{(m)}(a,\ell a)=0\) for \(a>R:=\|A\|+\sup_m\|B_m\|\) and
\(\ell\in U\). Both properties pass to the limit. Since \(K\) is
arbitrary and the lines \(b=\lambda_sa\) have measure zero, this proves
the statement on \(\mu_c\).

Let \(\lambda_{\min}(A)\le a\le R\), \(\ell\in K\) and \(m\ge m_0\). The
first factor of \(C_1^{(m)}(a,\ell a)\),
\(I-(A+\ell B_m)/(a(1+\ell^2))\), is bounded in norm by
\(1+R/\lambda_{\min}(A)\), and the displayed identity bounds the second
by a constant times \(1/a\). Hence \(\|C_1^{(m)}(a,\ell a)\|\le C/a\),
with \(C\) independent of \(m\), \(\ell\) and \(a\). Hence
\(\mu_c^{(m)}(a,\ell a)\le nC/(2\pi a)\) and \[
|g_m(\ell)|\le p\,|p-1|\,\frac{nC}{2\pi}\int_{\lambda_{\min}(A)}^{R}a^{-p}\,da ,
\] which is the uniform bound. The integrands vanish outside
\([\lambda_{\min}(A),R]\) and are dominated there by
\(\tfrac{nC}{2\pi}a^{-p}\), so dominated convergence gives
\(g_m(\ell)\to g(\ell)\).
\end{proof}

\subsection{Proof of the representation}\label{sec-proof-rep}

\begin{proof}
If \(B=0\), then \(L=0\) and the assertion is immediate with
\(\rho_p=\operatorname{tr}(A^{-p})\delta_0\). Henceforth assume
\(B\ne0\), so \(l_n>0\).

\emph{Step 1: uniqueness.} A finite signed measure on a bounded interval
is determined by its moments, and if \(\rho\) represents \(F_p\) then
\(F_p^{(k)}(0)=(-1)^k(p)_k\int\ell^k\,d\rho\) with
\((p)_k=p(p+1)\cdots(p+k-1)\neq0\). So \(\rho\) is determined by
\(F_p\). The same holds without a bound on the support. This is the
classical uniqueness of the representing measure of a generalized
Stieltjes transform (it is used, for instance, in
(\citeproc{ref-karp2012generalized}{Karp and Prilepkina 2012})). We
recall the argument. Let \(\sigma\) be a finite signed measure on
\([0,\infty)\) and \(G(t)=\int(1+t\ell)^{-p}\,d\sigma(\ell)\), \(t>0\).
Writing
\((1+t\ell)^{-p}=\Gamma(p)^{-1}\int_0^\infty u^{p-1}e^{-u(1+t\ell)}\,du\),
exchanging the integrals (the double integral of the absolute value is
at most \(\Gamma(p)\,|\sigma|([0,\infty))\)), and substituting
\(u=v/t\), \[
\Gamma(p)\,t^{p}\,G(t)=\int_0^\infty e^{-v/t}\,v^{p-1}\,\widehat\sigma(v)\,dv,\qquad
\widehat\sigma(v):=\int e^{-v\ell}\,d\sigma(\ell).
\] Thus \(G\) determines the Laplace transform of the locally integrable
function \(v^{p-1}\widehat\sigma(v)\), hence this function almost
everywhere, hence the continuous function \(\widehat\sigma\) everywhere.
Finally, \(\widehat\sigma\) is the Laplace transform of the finite
signed measure \(\sigma\), which it determines.

\emph{Step 2: simple eigenvalues.} Assume that the eigenvalues of \(L\)
are simple. For \(p\ne1\) and \(\varepsilon>0\) set \[
f_\varepsilon(x):=p(p-1)x^{-p}\mathbf 1_{\{x>\varepsilon\}}.
\] A direct integration gives, for \(x>\varepsilon\),
\begin{equation}\phantomsection\label{eq-truncated-inverse}{
H(f_\varepsilon)(x)
=x^{-p}+(p-1)\varepsilon^{-p}-p\varepsilon^{1-p}x^{-1}.
}\end{equation} Choose \(0<\varepsilon<\lambda_{\min}(A)\). Since
\(A+tB\succeq A\), every eigenvalue of \(X_t:=A+tB\) exceeds
\(\varepsilon\), and therefore
\begin{equation}\phantomsection\label{eq-trunc-left}{
\operatorname{tr}H(f_\varepsilon)(X_t)
=F_p(t)+n(p-1)\varepsilon^{-p}-p\varepsilon^{1-p}F_1(t).
}\end{equation} By (\ref{eq-heinavaara}) the left side is
\(\int f_\varepsilon(a+tb)\,d\mu_{A,B}\). Write
\(\mu_{A,B}=\mu_c+\sum_i\gamma_i\) as in (\ref{eq-singpart}).

The continuous part contributes
\begin{equation}\phantomsection\label{eq-trunc-cont}{
p(p-1)\int_{\{a+tb>\varepsilon\}}(a+tb)^{-p}\,d\mu_c(a,b),
}\end{equation} a finite integral: on \(\operatorname{supp}\mu_c\) one
has \(a,b\ge0\) by Lemma~\ref{lem-cone}(ii), so the domain of
integration is bounded away from the origin, where \(\mu_{A,B}\) is
finite by Lemma~\ref{lem-cone}(i). On the singular ray of slope \(l_i\),
(\ref{eq-singpart}) gives, with \(w_i:=1+tl_i\) and for every bounded
Borel \(f\) vanishing near the origin (to which (\ref{eq-singpart})
extends from smooth test functions, both sides being positive Radon
measures on \(\mathbb R^2\setminus\{0\}\)),
\begin{equation}\phantomsection\label{eq-ray-is-H}{
\int f(a+tb)\,d\gamma_i(a,b)
=\int_0^1\frac{1-s}{s}\,f\!\Bigl(\frac{s\,w_i}{N_i}\Bigr)ds
=H(f)\Bigl(\frac{w_i}{N_i}\Bigr).
}\end{equation} Since
\(w_i/N_i=\langle v_i|A+tB|v_i\rangle\ge\lambda_{\min}(A)>\varepsilon\)
for the normalized generalized eigenvector \(v_i\) of
(\ref{eq-singpart}), (\ref{eq-truncated-inverse}) applies at
\(x=w_i/N_i\) and gives
\begin{equation}\phantomsection\label{eq-trunc-ray}{
\int f_\varepsilon(a+tb)\,d\gamma_i
=N_i^p w_i^{-p}+(p-1)\varepsilon^{-p}-p\varepsilon^{1-p}\frac{N_i}{w_i}.
}\end{equation}

Summing (\ref{eq-trunc-ray}) over \(i\) and using the exact \(p=1\)
identity \[
F_1(t)=\sum_i\frac{N_i}{1+tl_i}
\] shows that both cutoff terms in (\ref{eq-trunc-left}) cancel. Hence,
for \(p\ne1\) and every \(0<\varepsilon<\lambda_{\min}(A)\),
\begin{equation}\phantomsection\label{eq-trunc-eps}{
F_p(t)=\sum_i\frac{N_i^p}{(1+tl_i)^p}
+p(p-1)\int_{\{a+tb>\varepsilon\}}(a+tb)^{-p}\,d\mu_c(a,b).
}\end{equation} The left side and the sum do not depend on
\(\varepsilon\), so neither does the integral. Its integrand is
positive, so \(\mu_c\) gives no mass to
\(\{\varepsilon_1<a+tb\le\varepsilon_2\}\) for
\(0<\varepsilon_1<\varepsilon_2<\lambda_{\min}(A)\), and hence none to
\(\{0<a+tb<\lambda_{\min}(A)\}\). The cone (\ref{eq-cone}) meets
\(\{a+tb\le0\}\) only at the origin, so the indicator in
(\ref{eq-trunc-eps}) can be dropped:
\begin{equation}\phantomsection\label{eq-simple-trunc-rep}{
F_p(t)=\sum_i\frac{N_i^p}{(1+tl_i)^p}
+p(p-1)\int(a+tb)^{-p}\,d\mu_c(a,b).
}\end{equation} At \(t=0\) the same observation gives
\begin{equation}\phantomsection\label{eq-support-bound}{
a\ge\lambda_{\min}(A)\qquad\text{for }\mu_c\text{-almost every }(a,b).
}\end{equation} This is a statement about \(\mu_c\) alone: it has been
obtained with an auxiliary exponent \(p\neq1\) and does not depend on
the exponent of the theorem. For \(p=1\) the identity for \(F_1\)
displayed above is the representation. Using \(b=\ell a\) on the
continuous part and (\ref{eq-cone}), \[
p(p-1)\int(a+tb)^{-p}\,d\mu_c
=\int_{l_1}^{l_n}(1+t\ell)^{-p}
 p(p-1)\Bigl[\int_0^\infty a^{1-p}\mu_c(a,\ell a)\,da\Bigr]d\ell,
\] which proves (\ref{eq-main-rep}) and (\ref{eq-cts}) for simple
eigenvalues. By (\ref{eq-support-bound}) the radial integral in
(\ref{eq-cts}) extends only over \(a\ge\lambda_{\min}(A)\), and it
converges.

\emph{Step 3: repeated eigenvalues.} Choose positive semidefinite
matrices \(L_m\to L\) with simple spectrum and set
\(B_m:=A^{1/2}L_mA^{1/2}\to B\). By Step 2, \[
F_p^{(m)}(t)=\sum_i\frac{(N_i^{(m)})^{p}}{(1+tl_i^{(m)})^{p}}+\int(1+t\ell)^{-p}g_m(\ell)\,d\ell ,
\] with \(g_m\) the density (\ref{eq-cts}) of \((A,B_m)\), and
\(\mu_c^{(m)}\) satisfies (\ref{eq-support-bound}). As \(m\to\infty\),
\(F_p^{(m)}\to F_p\) pointwise. The atomic parts have weights bounded by
\(\|A^{-1}\|^p\) and nodes \(l_i^{(m)}\to\lambda_s\) for \(i\in I_s\),
the index set of the \(k_s\) eigenvalues converging to \(\lambda_s\), so
along a subsequence they converge weakly to
\(\sum_sc_s\delta_{\lambda_s}\). By Lemma~\ref{lem-density-limit},
\(\int\varphi\,g_m\to\int\varphi\,g\) for every continuous \(\varphi\)
supported away from \(\{\lambda_s\}\), where \(g\) is the density
(\ref{eq-cts}) of \((A,B)\), whose radial integral again extends only
over \(a\ge\lambda_{\min}(A)\). The total variations of the
\(g_m\,d\ell\) are bounded (their masses are
\(\operatorname{tr}A^{-p}-\sum_i(N_i^{(m)})^{p}\), of one sign). Hence
\(g\) is integrable, by Fatou's lemma, and along a further subsequence
\(g_m\,d\ell\to g\,d\ell+\varsigma\) weakly with \(\varsigma\) carried
by the finite set \(\{\lambda_s\}\). The limit
\(\rho_p:=\sum_s(c_s+\varsigma(\{\lambda_s\}))\delta_{\lambda_s}+g\,d\ell\)
represents \(F_p\), since \((1+t\ell)^{-p}\) is continuous on the common
compact support. By Step 1 it does not depend on the subsequence. This
is (\ref{eq-main-rep}) with \(\rho_p^{\mathrm{ac}}=g\,d\ell\) and atomic
weights \(\omega_s\) still to be identified.

\emph{Step 4: the atomic weights.} Both sides of (\ref{eq-main-rep}) are
analytic in \(t\) on an open connected neighbourhood of
\(\{\operatorname{Im}t>0\}\cup(-1/l_n,\infty)\), where the spectrum of
\(A+tB\) and the numbers \(1+t\ell\) avoid \((-\infty,0]\), so the
identity holds there. For \(\operatorname{Im}t>0\) put \(z:=-1/t\),
\(\operatorname{Im}z>0\): since \(A+tB=(-t)(zA-B)\) and
\(1+t\ell=(-t)(z-\ell)\), with all spectra in the open upper half-plane
and \(\arg(-t)\in(-\pi,0)\), the principal powers factor,
\((A+tB)^{-p}=z^{p}(zA-B)^{-p}\) and
\((1+t\ell)^{-p}=z^{p}(z-\ell)^{-p}\), and (\ref{eq-main-rep}), divided
by \(z^p\), becomes the identity (\ref{eq-stieltjes-p}) for
\(\Phi_p(z)=\operatorname{tr}[(zA-B)^{-p}]\), with
\(\eta=\rho_p^{\mathrm{ac}}\), finite and absolutely continuous. By
Lemma~\ref{lem-Apole}(ii),
\(\omega_s=\operatorname{tr}(\widetilde N_s^{\,p})=\Omega_s\) for every
\(p>0\): this is (\ref{eq-atomic-Omega}). For a simple eigenvalue it
recovers \(\Omega_s=N_s^{p}\), the weight found above from the singular
ray of (\citeproc{ref-heinavaara2024tracial}{Heinävaara 2025}, Thm.
1.6).

\emph{Step 5: masses.} Setting \(t=0\) gives the total mass
\(F_p(0)=\operatorname{tr}(A^{-p})\). Matching the coefficient of \(t\)
against \((1+t\ell)^{-p}=1-pt\ell+O(t^2)\) and
\(F_p(t)=\operatorname{tr}(A^{-p})-pt\operatorname{tr}(A^{-p-1}B)+O(t^2)\)
gives the first moment. The atomic mass is
\(\sum_s\Omega_s=\operatorname{tr}[(\mathcal C_LN)^p]\)
(block-diagonality of \(\mathcal C_LN\)), and subtracting yields
(\ref{eq-mass-identity}). Its sign is \(\operatorname{sgn}(p-1)\) by
Lemma~\ref{lem-pinching}.

\emph{Step 6: positivity.} For \(p>1\), \(p(p-1)>0\) and \(\mu_c\ge0\)
give \(\rho_p^{\mathrm{ac}}\ge0\), so \(\rho_p\ge0\). For \(p=1\) the
prefactor vanishes and \(\rho_p=\rho_p^{\mathrm{at}}\ge0\). For
\(0<p<1\), \(p(p-1)<0\) gives \(\rho_p^{\mathrm{ac}}\le0\), which is
absolutely continuous and so equals \(\rho_p\) off the atoms. Hence
\(\rho_p\ge0\) iff \(\rho_p^{\mathrm{ac}}=0\) iff its mass
(\ref{eq-mass-identity}) vanishes, iff (by Lemma~\ref{lem-pinching},
\(p\neq1\)) \([A,B]=0\). When \([A,B]=0\), \(A,B\) are jointly
diagonalizable and (\ref{eq-main-rep}) is the commuting atomic sum
\(\sum_ja_j^{-p}\delta_{b_j/a_j}\). A positive measure on \([0,\infty)\)
representing \(F_p\) is finite, with mass \(F_p(0)\), so by Step 1 it is
\(\rho_p\). Hence, for \(0<p<1\) and \([A,B]\neq0\), the function
\(F_p\) is not a Stieltjes function of order \(p\).
\end{proof}

\begin{corollary}[Integer \(p\): the two forms of the
density]\protect\hypertarget{cor-bspline-bridge}{}\label{cor-bspline-bridge}

For every integer \(p\ge2\) and almost every \(\ell\),
\begin{equation}\phantomsection\label{eq-bspline-bridge}{
p(p-1)\int_0^\infty a^{1-p}\mu_c(a,\ell a)\,da
=\!\!\sum_{\substack{i_1,\dots,i_p\\ l_{i_1},\dots,l_{i_p}\text{ not all equal}}}\!\! N_{i_1i_2}\cdots N_{i_pi_1}\,M(\ell\mid l_{i_1},\dots,l_{i_p})
\qquad(p\in\mathbb Z_{\ge2}).
}\end{equation}

\end{corollary}

\begin{proof}
By uniqueness, the measure of Theorem~\ref{thm-rep} is the B-spline
measure of Theorem~\ref{thm-bspline}. The two sides are the densities
(\ref{eq-cts}) and (\ref{eq-bspline-density}) of its continuous part.
\end{proof}

For integer \(p\) the radial integral of Heinävaara's density is thus a
piecewise polynomial in \(\ell\), with coefficients given by the
eigendata of \((A,B)\), computed with no analysis of the pencil
spectrum.

\begin{refremark}[A second proof, by a gamma integral]
For simple eigenvalues the representation (\ref{eq-simple-trunc-rep})
can also be obtained from the gamma integral
\(x^{-p}=\Gamma(p+1)^{-1}\int_0^\infty\tau^p\,xe^{-\tau x}\,d\tau\), a
variant of the device by which Lieb and Seiringer pass from the Laplace
transform to the functions \(F_p\)
(\citeproc{ref-lieb-seiringer2004}{Lieb and Seiringer 2004}, proof of
Thm. 1). It gives
\(F_p(t)=\Gamma(p+1)^{-1}\int_0^\infty\tau^p\operatorname{tr}[X_te^{-\tau X_t}]\,d\tau\),
and the function \(x\mapsto xe^{-\tau x}\) has an explicit admissible
preimage under \(H\): from \(H(y^{n+1})=y^{n+1}/((n+1)(n+2))\) one finds
\(H\bigl(y\,\psi(\tau y)\bigr)(x)=xe^{-\tau x}\) with
\(\psi(v)=(v^2-4v+2)e^{-v}\). By (\ref{eq-heinavaara}), and by
(\ref{eq-ray-is-H}), which holds for this function as well (by dominated
convergence), one has for every \(\tau>0\), with \(y=a+tb\), \[
\operatorname{tr}\bigl[X_te^{-\tau X_t}\bigr]
=\int y\,\psi(\tau y)\,d\mu_c(a,b)+\sum_i\frac{w_i}{N_i}\,e^{-\tau w_i/N_i}.
\] Integrate against \(\Gamma(p+1)^{-1}\tau^p\,d\tau\). The sum gives
\(\sum_iN_i^pw_i^{-p}\). On the continuous part, Fubini's theorem and
\(\Gamma(p+1)^{-1}\int_0^\infty\tau^p\,y\psi(\tau y)\,d\tau=p(p-1)\,y^{-p}\)
give \(p(p-1)\int(a+tb)^{-p}\,d\mu_c\), which is
(\ref{eq-simple-trunc-rep}). The exchange of integrals requires
\(\int(a+tb)^{-p}\,d\mu_c<\infty\), for which a lower bound for \(a\) on
\(\operatorname{supp}\mu_c\) suffices. Such a bound is Heinävaara's
support statement (\citeproc{ref-heinavaara-thesis}{Heinävaara 2024b},
Prop. 4). The cutoff argument of Step 2 does not need it, and yields it
(\ref{eq-support-bound}). In the present argument the coefficient
\(p(p-1)\) appears as
\([\Gamma(p+3)-4\Gamma(p+2)+2\Gamma(p+1)]/\Gamma(p+1)\).

\label{rem-gamma-route}

\end{refremark}

\section{\texorpdfstring{The \(n=2\) case in closed
form}{The n=2 case in closed form}}\label{sec-n2}

For \(n=2\) every object collapses to elementary closed form and the
whole family \(\operatorname{tr}(A+tB)^{-p}\), \(0<p\le\infty\), can be
written out. The auxiliary matrix (\ref{eq-M-entries}) is \(2\times2\),
so the imaginary spectrum driving the continuous density is a single
square root: the finite-\(p\) density is a Gauss hypergeometric
\({}_2F_1\), its integer-\(p\) specializations are the B-splines of
Section~\ref{sec-bsplines}, and its \(p\to\infty\) confluence is the
Bessel-\(I_1\) density of Stahl's \(2\times2\) measure
(Section~\ref{sec-n2-stahl}). One complex number
\(N_{12}=\langle e_1|A^{-1}|e_2\rangle\) measures non-commutativity
throughout, vanishing iff \([A,B]=0\).

\subsection{Ray coordinates and the auxiliary
matrix}\label{sec-auxiliary}

In ray coordinates \(b=\ell a\) the matrix \(C_1(a,b)\) of
(\ref{eq-rho-c}), whose imaginary spectrum carries \(\mu_c\), is similar
to an affine function of \(a\). This is the tool for the closed forms
below and for Section~\ref{sec-stahl}.

\begin{lemma}[The auxiliary
matrix]\protect\hypertarget{lem-MatrxM}{}\label{lem-MatrxM}

With \(b=\ell a\), set \(K(a,\ell):=a^{2}A^{-1/2}C_1(a,\ell a)A^{1/2}\),
\(C_1\) as in (\ref{eq-rho-c}). In the \(L\)-eigenbasis \(\{e_i\}\) its
entries are \begin{equation}\phantomsection\label{eq-M-entries}{
K_{jk}(a,\ell)=d_j(\ell)\,\delta_{jk}-a\,\frac{N_{jk}}{l_k-\ell},
\qquad
d_j(\ell):=\frac{1+\ell l_j}{(1+\ell^2)(l_j-\ell)},\qquad N_{jk}=\langle e_j|A^{-1}|e_k\rangle .
}\end{equation} For fixed \(\ell\) the map
\(l\mapsto d(l;\ell):=\frac{1+\ell l}{(1+\ell^2)(l-\ell)}\) is a Möbius
transformation, injective on \(l\neq\ell\). Hence the diagonal entries
\(d_j(\ell)=d(l_j;\ell)\) are pairwise distinct whenever the \(l_j\)
are.

\end{lemma}

\begin{proof}
In ray coordinates \(b A-a B=a A^{1/2}(\ell I-L)A^{1/2}\) and
\(a A+b B=a A^{1/2}(I+\ell L)A^{1/2}\), with \(a^2+b^2=a^2(1+\ell^2)\).
Hence \[
a^{2}A^{-1/2}C_1(a,\ell a)A^{1/2}
=a A^{-1}(\ell I-L)^{-1}-\frac{I+\ell L}{1+\ell^2}(\ell I-L)^{-1}.
\] Since
\((\ell I-L)^{-1}=\operatorname{diag}\big((\ell-l_k)^{-1}\big)\) in the
\(L\)-eigenbasis, the second term is diagonal with entries
\(-\tfrac{1+\ell l_j}{(1+\ell^2)(\ell-l_j)}=d_j(\ell)\) and the first
has entries \(a N_{jk}(\ell-l_k)^{-1}=-a N_{jk}/(l_k-\ell)\), which is
(\ref{eq-M-entries}). The map \(d(\,\cdot\,;\ell)\) is the Möbius
transformation with matrix
\(\left(\begin{smallmatrix}\ell & 1\\ 1+\ell^2 & -\ell(1+\ell^2)\end{smallmatrix}\right)\),
of determinant \(-(1+\ell^2)^2\neq0\), hence injective on \(l\neq\ell\).
\end{proof}

Since \(K(a,\ell)\) is similar to \(a^{2}C_1(a,\ell a)\), its
eigenvalues are \(a^{2}\) times those of \(C_1(a,\ell a)\), so in ray
coordinates \begin{equation}\phantomsection\label{eq-mu-ray}{
\mu_c(a,\ell a)=\frac{1}{2\pi a^{2}}\sum_{k=1}^n\bigl|\operatorname{Im}\lambda_k[K(a,\ell)]\bigr| ,
}\end{equation} and the continuous density is supported where \(K\) has
non-real eigenvalues.

\subsection{\texorpdfstring{Finite \(p>1\): the Gauss hypergeometric
density}{Finite p\textgreater1: the Gauss hypergeometric density}}\label{finite-p1-the-gauss-hypergeometric-density}

Work in the \(L\)-eigenbasis, \(L=\operatorname{diag}(l_1,l_2)\),
\(A^{-1}=\bigl(\begin{smallmatrix}N_1&N_{12}\\\bar N_{12}&N_2\end{smallmatrix}\bigr)\)
with \(N_1N_2-|N_{12}|^2=\det A^{-1}>0\). With \(d_j(\ell)\) as in
(\ref{eq-M-entries}) put, for \(\ell\in(l_1,l_2)\),
\begin{equation}\phantomsection\label{eq-n2p-PQK}{
P=d_1(\ell)-d_2(\ell),\quad Q=\frac{N_1}{l_1-\ell}-\frac{N_2}{l_2-\ell},\quad
\kappa=\frac{|N_{12}|^2}{(l_1-\ell)(l_2-\ell)}<0,\quad D=Q^2-4|\kappa| .
}\end{equation} Writing \(x:=\ell-l_1>0\) and \(y:=l_2-\ell>0\), one has
\(Q=-(N_1/x+N_2/y)<0\) and, since
\((1+\ell l_1)y+(1+\ell l_2)x=(l_2-l_1)(1+\ell^2)\),
\(P=-(l_2-l_1)/(xy)<0\). So \(PQ>0\) on all of \((l_1,l_2)\), and
\(D=Q^2-4|\kappa|>0\) there (see the proof below).

\begin{proposition}[Finite-\(p\) density,
\(n=2\)]\protect\hypertarget{prp-n2-finitep}{}\label{prp-n2-finitep}

For \(n=2\), \(A\succ0\), \(B\succeq0\) with \(l_1<l_2\) and \(p>1\),
\(\rho_p=N_1^{p}\delta_{l_1}+N_2^{p}\delta_{l_2}+g_p(\ell)\,d\ell\) with
continuous density on \((l_1,l_2)\) given, in terms of \(P,Q,\kappa,D\)
of (\ref{eq-n2p-PQK}), by
\begin{equation}\phantomsection\label{eq-n2p-density}{
g_p(\ell)=\frac{p(p-1)}{4}\sqrt{D}\,\varrho^2\,\bar a^{-1-p}\,
{}_2F_1\!\Bigl(\tfrac{1+p}{2},\tfrac{2+p}{2};2;\tfrac{4|\kappa|}{Q^2}\Bigr),
\qquad \bar a=\frac{PQ}{D},\ \ \varrho^2=\frac{4P^2|\kappa|}{D^2},
}\end{equation} where \(4|\kappa|/Q^2\in[0,1)\) and
\({}_2F_1(a,b;c;z)=\sum_{k\ge0}\frac{(a)_k(b)_k}{(c)_k}\frac{z^k}{k!}\)
is the Gauss hypergeometric series, \((a)_k=a(a+1)\cdots(a+k-1)\).

\end{proposition}

\begin{proof}
The matrix
\(K(a,\ell)=\operatorname{diag}(d_1,d_2)-a A^{-1}(L-\ell I)^{-1}\) of
(\ref{eq-M-entries}) is affine in \(a\). Its eigenvalues are real or a
conjugate pair, and in the latter case
\(\sum_i|\operatorname{Im}\lambda_i[K]|=\sqrt{-\Delta_K(a)}\) with
\(\Delta_K(a)=(\operatorname{tr}K)^2-4\det K=P^2-2PQ\,a+D\,a^2\). Here
\(D=(\operatorname{tr}K^{(1)})^2-4\det K^{(1)}\),
\(K^{(1)}=A^{-1}(L-\ell I)^{-1}\), and since
\(\det K^{(1)}=\det(A^{-1})/[(l_1-\ell)(l_2-\ell)]<0\) on \((l_1,l_2)\),
\(D>0\). The upward parabola \(\Delta_K\) has \(\Delta_K(0)=P^2>0\) and,
as \(PQ>0\), dips below zero exactly on \((a_-,a_+)\),
\(a_\pm=(PQ\pm2|P|\sqrt{|\kappa|})/D\), with \(0<a_-<a_+\). By
(\ref{eq-cts}) and
\(|\operatorname{Im}\lambda_i[C_1]|=a^{-2}|\operatorname{Im}\lambda_i[K]|\),
\[
g_p(\ell)=\frac{p(p-1)}{2\pi}\sqrt{D}\int_{a_-}^{a_+}a^{-1-p}\sqrt{(a-a_-)(a_+-a)}\,da .
\] The substitution \(a=\bar a+\varrho\cos\theta\)
(\(\bar a=(a_-+a_+)/2=PQ/D\),
\(\varrho=(a_+-a_-)/2=2|P|\sqrt{|\kappa|}/D\)) gives a Gegenbauer
integral, \[
\int_{a_-}^{a_+}a^{-1-p}\sqrt{(a-a_-)(a_+-a)}\,da
=\varrho^2\!\int_0^\pi\frac{\sin^2\theta\,d\theta}{(\bar a+\varrho\cos\theta)^{1+p}}
=\frac{\pi}{2}\varrho^2\bar a^{-1-p}\,{}_2F_1\!\Bigl(\tfrac{1+p}{2},\tfrac{2+p}{2};2;\tfrac{\varrho^2}{\bar a^2}\Bigr)
\] valid for \(\varrho/\bar a<1\): expand
\((\bar a+\varrho\cos\theta)^{-1-p}\) in powers of
\(\varrho\cos\theta/\bar a\) and use
\(\int_0^\pi\sin^2\theta\cos^{2k}\theta\,d\theta=\tfrac\pi2\,(\tfrac12)_k/(2)_k\),
the odd powers integrating to zero (equivalently,
(\citeproc{ref-gradshteyn-ryzhik2015}{Gradshteyn and Ryzhik 2015,
3.665.2}) followed by a quadratic transformation of the hypergeometric
function). Since \(\varrho^2/\bar a^2=4|\kappa|/Q^2\) and
\(D=Q^2-4|\kappa|>0\), this holds and yields (\ref{eq-n2p-density}).
\end{proof}

Every factor is nonnegative for \(p>1\), so the density is manifestly
positive, as Heinävaara's density is at every order, including
\(p=\infty\) (Section~\ref{sec-stahl}).

\textbf{Integer \(p\): reduction to B-splines.} For integer \(p\) the
series does not terminate but sums to an algebraic function ---
\({}_2F_1(\tfrac32,2;2;z)=(1-z)^{-3/2}\) at \(p=2\),
\({}_2F_1(2,\tfrac52;2;z)=(1-z)^{-5/2}\) at \(p=3\) --- and the square
roots cancel against the algebraic prefactors: with
\(z=\varrho^2/\bar a^2\) one has
\(\bar a^2(1-z)=\bar a^2-\varrho^2=a_-a_+=P^2/D\), so at \(p=2\) \[
g_2(\ell)=\tfrac12\sqrt D\,\varrho^2\,(P^2/D)^{-3/2}=\frac{2|\kappa|}{|P|}=\frac{2|N_{12}|^2}{l_2-l_1},
\] constant on \((l_1,l_2)\): the box
\(2|N_{12}|^2M(\cdot\mid l_1,l_2)\) of the two words \((1,2)\),
\((2,1)\) in (\ref{eq-bspline-density}), the \(n=2\) instance of the
bridge (\ref{eq-bspline-bridge}). In the same way (\ref{eq-n2p-density})
is a polynomial of degree \(p-2\) on \((l_1,l_2)\) for every integer
\(p\), as Theorem~\ref{thm-bspline} requires.

\section{\texorpdfstring{The limit \(p\to\infty\): Stahl's measure as a
confluence}{The limit p\textbackslash to\textbackslash infty: Stahl's measure as a confluence}}\label{sec-stahl}

The representation of Theorem~\ref{thm-rep} holds for every \(p>0\).
Letting \(p\to\infty\) along the \emph{Hille--Yosida ray}
\(A\rightsquigarrow A_p:=pI+A_0\) (\(A_0=A_0^*\), \(B\succeq0\) fixed)
exhibits Stahl's measure of \(\operatorname{tr}\,e^{-(A_0+tB)}\) as the
limit of the finite-\(p\) measures \(\rho_p\equiv\rho_{A_p,B,p}\). Its
continuous density is the manifestly positive real form of Heinävaara's
thesis. We recover it here as a confluence rather than reprove it.

The passage from transforms of order \(p\) to the Laplace transform,
\((1+tb/p)^{-p}\to e^{-tb}\), is the device by which Curry and
Schoenberg determined the limits of B-splines of growing order
(\citeproc{ref-curry-schoenberg1966}{Curry and Schoenberg 1966, pt.
II}). Here the same limit is taken for the measures \(\rho_p\), which
for integer \(p\) are sums of B-splines.

Throughout this section \(B\) has simple spectrum
\(0\le b_1<\cdots<b_n\) with eigenvectors \(u_i\), and
\(a_i:=\langle u_i|A_0|u_i\rangle\). For \(p>\|A_0\|\),
Theorem~\ref{thm-rep} applies to \((A_p,B)\). Write
\(L_p,e_i^{(p)},N_i^{(p)}\) for its spectral data and define the
\textbf{rescaled measure}
\begin{equation}\phantomsection\label{eq-mup-def}{
\widetilde\rho_p:=p^{\,p}(S_p)_*\rho_p,\qquad S_p(\ell):=p\ell ,
}\end{equation} supported on \([p\,l_1^{(p)},p\,l_n^{(p)}]\).

\begin{theorem}[Stahl
limit]\protect\hypertarget{thm-stahl-limit}{}\label{thm-stahl-limit}

Let \(A_0=A_0^*\), let \(B\succeq0\) have simple spectrum
\(0\le b_1<\dots<b_n\) with eigenvectors \(u_i\),
\(a_i:=\langle u_i|A_0|u_i\rangle\), and let \(\widetilde\rho_p\) be the
rescaled measures (\ref{eq-mup-def}) of the pairs \((pI+A_0,B)\),
\(p>\|A_0\|\). Let \(\rho_\infty\) denote Stahl's measure, the unique
finite measure on \([b_1,b_n]\) with
\begin{equation}\phantomsection\label{eq-stahl-laplace}{
\operatorname{tr}\,e^{-(A_0+tB)}=\int e^{-tb}\,d\rho_\infty(b)\qquad(t\ge0).
}\end{equation}

\begin{enumerate}
\def\labelenumi{(\roman{enumi})}
\item
  \(\widetilde\rho_p\to\rho_\infty\) weakly as \(p\to\infty\). In
  particular, \(\rho_\infty\ge0\).
\item
  For every \(b\notin\{b_1,\dots,b_n\}\) the density of
  \(\widetilde\rho_p^{\mathrm{ac}}\) converges,
  \(\widetilde\rho_p^{\mathrm{ac}}(b)\to\rho_\infty^{\mathrm{ac}}(b)\),
  where \begin{equation}\phantomsection\label{eq-phi-def}{
  \rho_\infty^{\mathrm{ac}}(b)=\frac{1}{2\pi}\int_{\mathbb R}e^{\tau}\sum_{i=1}^n\bigl|\operatorname{Im}\lambda_i[(A_0+\tau I)(B-bI)^{-1}]\bigr|\,d\tau ,
  }\end{equation} a manifestly nonnegative function supported in
  \([b_1,b_n]\).
\item
  \(\rho_\infty=\sum_{i=1}^n e^{-a_i}\delta_{b_i}+\rho_\infty^{\mathrm{ac}}(b)\,db\).
\end{enumerate}

\end{theorem}

The density (\ref{eq-phi-def}) is Heinävaara's thesis form of Stahl's
measure (\citeproc{ref-heinavaara-thesis}{Heinävaara 2024b, sec. 1.4},
proof of Thm. 24), and coincides by Laplace uniqueness with Stahl's
contour representation (\citeproc{ref-stahl2013}{Stahl 2013};
\citeproc{ref-eremenko2015}{Eremenko 2015}, Prop. 1). The atoms in (iii)
are in both sources. What Theorem~\ref{thm-stahl-limit} adds is the
convergence: measure and density arise as limits of the finite-\(p\)
objects, so that the positivity of \(\rho_\infty\) appears as a limit of
that of the \(\rho_p\). That Stahl's theorem follows from positivity at
finite \(p\) is one of the equivalences of Lieb and Seiringer
(\citeproc{ref-lieb-seiringer2004}{Lieb and Seiringer 2004}, Thm. 1),
which they prove through the polynomial form of the statement. Since
Theorem~\ref{thm-rep} itself rests on Heinävaara's theorem, the limit is
not an independent proof of Stahl's theorem.

\begin{lemma}[Trace limit]\protect\hypertarget{lem-HY}{}\label{lem-HY}

For each \(t\ge0\),
\(\int(1+tb/p)^{-p}\,d\widetilde\rho_p(b)=p^{\,p}\operatorname{tr}(A_p+tB)^{-p}\to\operatorname{tr}\,e^{-(A_0+tB)}\),
and \(\widetilde\rho_p(\mathbb R)\to\operatorname{tr}\,e^{-A_0}\).

\end{lemma}

\begin{proof}
The equality is (\ref{eq-mup-def}) with Theorem~\ref{thm-rep} at
\((A_p,B)\). Writing \(A_p+tB=p(I+X_t/p)\), \(X_t:=A_0+tB\), and
diagonalizing reduces the limit to \((1+x/p)^{-p}\to e^{-x}\) on the
eigenvalues of \(X_t\). The mass is the case \(t=0\).
\end{proof}

\begin{lemma}[Continuous
limit]\protect\hypertarget{lem-cts}{}\label{lem-cts}

In the setting of Theorem~\ref{thm-stahl-limit}, fix
\(b\notin\{b_1,\dots,b_n\}\). Then
\(\widetilde\rho_p^{\mathrm{ac}}(b)\to\rho_\infty^{\mathrm{ac}}(b)\)
with \(\rho_\infty^{\mathrm{ac}}\) as in (\ref{eq-phi-def}), and
\(\rho_\infty^{\mathrm{ac}}(b)=0\) for \(b\notin[b_1,b_n]\). Moreover,
for every compact \(J'\subset\mathbb R\setminus\{b_1,\dots,b_n\}\) the
densities \(\widetilde\rho_p^{\mathrm{ac}}\) are bounded on \(J'\),
uniformly in \(p\) large.

\end{lemma}

\begin{proof}
\emph{Reduction to a \(\tau\)-integral.} By (\ref{eq-cts}),
(\ref{eq-mu-ray}) and the Jacobian \(d\ell=db/p\) of \(S_p\),
\(\widetilde\rho_p^{\mathrm{ac}}(b)=p^{\,p-1}\rho_p^{\mathrm{ac}}(b/p)
=\frac{p^{\,p}(p-1)}{2\pi}\int_0^\infty a^{-1-p}\,\Sigma_p^{\circ}(a,b)\,da\)
with
\(\Sigma_p^{\circ}(a,b):=\sum_i|\operatorname{Im}\lambda_i[K_p(a,b/p)]|\),
where \(K_p\) is the matrix \(K\) of Lemma~\ref{lem-MatrxM} for the pair
\((A_p,B)\). Substituting \(a=p-\tau\) cancels every divergent power of
\(p\): \begin{equation}\phantomsection\label{eq-psip-tau}{
\widetilde\rho_p^{\mathrm{ac}}(b)=\frac{p-1}{2\pi p}\int_{-\infty}^{p}\Bigl(1-\frac{\tau}{p}\Bigr)^{-1-p}\Sigma_p(\tau,b)\,d\tau,
\qquad \Sigma_p(\tau,b):=\Sigma_p^{\circ}(p-\tau,b).
}\end{equation}

\emph{Pencil form.} By Lemma~\ref{lem-MatrxM}, \(K_p(a,\ell)\) with
\(a=p-\tau\), \(\ell=b/p\) is similar to
\(a^2C_1(a,\ell a)=(a^2Y_1)Y_2^{-1}\), where \(C_1\) is the matrix of
(\ref{eq-rho-c}) for the pair \((A_p,B)\) and, with \(\tilde b=\ell a\),
\(Y_1:=I-(aA_p+\tilde bB)/(a^2+\tilde b^2)\) and
\(Y_2:=\tilde bA_p-aB\). Using \(A_p=pI+A_0\), \(b=p\ell\) and
\(a^2/(a^2+\tilde b^2)=1/(1+\ell^2)\), one finds
\begin{equation}\phantomsection\label{eq-pencil}{
\begin{aligned}
Y_2&=-(p-\tau)\,\mathcal B_p,\qquad a^2Y_1=-(p-\tau)\,\mathcal A_p,\\
\mathcal A_p&:=(A_0+\tau I)+\frac{\ell(B-bI)-\ell^2A_0}{1+\ell^2},\qquad
\mathcal B_p:=(B-bI)-\ell A_0 ,
\end{aligned}
}\end{equation} the second identity resting on the cancellation
\(a-\tfrac{p}{1+\ell^2}=\tfrac{b\ell}{1+\ell^2}-\tau\) of the divergent
scalar. As \(p\to\infty\), \(\mathcal A_p\to A_0+\tau I\) and
\(\mathcal B_p\to B-bI\) in norm, uniformly for \((\tau,b)\) in compact
sets. Since \(b\notin\{b_1,\dots,b_n\}\), \(\mathcal B_p\) is invertible
for \(p\) large. Then \(K_p\) is similar to
\(\mathcal A_p\mathcal B_p^{-1}\), and \(\operatorname{spec}K_p\) is the
set of generalized eigenvalues of the Hermitian pencil
\((\mathcal A_p,\mathcal B_p)\).

\emph{Pointwise limit.} Since \(B-bI\) is invertible,
\(\mathcal B_p^{-1}\mathcal A_p\to(B-bI)^{-1}(A_0+\tau I)\) in norm,
which is similar to \(K_\infty(\tau,b):=(A_0+\tau I)(B-bI)^{-1}\). As
\(Z\mapsto\sum_i|\operatorname{Im}\lambda_i[Z]|\) is continuous,
\(\Sigma_p(\tau,b)\to\Sigma_\infty(\tau,b):=\sum_i|\operatorname{Im}\lambda_i[K_\infty(\tau,b)]|\),
and the integrand of (\ref{eq-psip-tau}) converges pointwise to
\(\tfrac1{2\pi}e^{\tau}\Sigma_\infty(\tau,b)\).

\emph{A reality window.} We claim \(\Sigma_p(\tau,b)=0\) for
\(\tau\notin[-\tau_-,\tau_+]\) once \(p\ge p_0\), with
\(\tau_-:=\|A_0\|+\|B\|\), \(\tau_+:=\|A_0\|+1\) and
\(p_0:=\max(\|A_0\|+\tau_+,\,1+|b|(|b|+\|B\|))\). The window does not
depend on \(b\), and \(p_0\) is bounded for \(b\) in a bounded set. The
lower cut-off is Lemma~\ref{lem-cone}(iii): \(\mu_c\) vanishes for
\(a>\|A_p\|+\|B\|\), and \(\|A_p\|+\|B\|\le p+\tau_-\), so it vanishes
for \(\tau<-\tau_-\). For the upper cut-off note that the generalized
eigenvalues of the Hermitian pencil (\ref{eq-pencil}) are real as soon
as \(\mathcal A_p\succ0\) (if \(\mathcal A_pw=\lambda\mathcal B_pw\)
with \(w\ne0\) then \(w^*\mathcal B_pw\neq0\) and
\(\lambda=w^*\mathcal A_pw/w^*\mathcal B_pw\in\mathbb R\)). Multiplying
(\ref{eq-pencil}) by \((1+\ell^2)(p-\tau)\) and writing
\(\tilde b=\ell(p-\tau)\), \[
(1+\ell^2)(p-\tau)\,\mathcal A_p=(p-\tau)(A_0+\tau I)-\tilde b(\tilde bI-B).
\] For \(\tau>\|A_0\|\) the first term is
\(\succeq(p-\tau)(\tau-\|A_0\|)I\), while
\(\|\tilde b(\tilde bI-B)\|\le\frac{p-\tau}{p}\,|b|(|b|+\|B\|)\) because
\(|\tilde b|=|b|(p-\tau)/p\le|b|\). Hence \(\mathcal A_p\succ0\)
whenever \(p(\tau-\|A_0\|)>|b|(|b|+\|B\|)\), in particular for all
\(\tau\ge\tau_+\) and \(p\ge p_0\).

\emph{Domination.} On the window \([-\tau_-,\tau_+]\) the convergence
\(\mathcal A_p\mathcal B_p^{-1}\to K_\infty(\tau,b)\) is uniform, so
\(\Sigma_p(\tau,b)\le n\|\mathcal A_p\mathcal B_p^{-1}\|\le C(b)\) for
\(p\) large, with \(C(b)\) bounded on compact subsets of
\(\mathbb R\setminus\{b_1,\dots,b_n\}\), and
\((1-\tau/p)^{-1-p}\le2e^{\tau}\) there for \(p\) large. The integrand
of (\ref{eq-psip-tau}) is therefore dominated by
\(\tfrac1\pi C(b)e^{\tau}\mathbf 1_{[-\tau_-,\tau_+]}(\tau)\), and
dominated convergence gives
\(\widetilde\rho_p^{\mathrm{ac}}(b)\to\frac1{2\pi}\int_{-\tau_-}^{\tau_+}e^{\tau}\Sigma_\infty(\tau,b)\,d\tau
=\rho_\infty^{\mathrm{ac}}(b)\), the window being also where
\(\Sigma_\infty\) lives (limit of the reality statement). The dominating
function is bounded uniformly for \(b\) in a compact subset of
\(\mathbb R\setminus\{b_1,\dots,b_n\}\), which gives the last assertion.

\emph{Support.} If \(b\notin[b_1,b_n]\) then \(B-bI\) is definite and,
with \(S:=|B-bI|^{1/2}\), \(K_\infty(\tau,b)\) is similar to the
Hermitian matrix \(\pm S^{-1}(A_0+\tau I)S^{-1}\). So
\(\Sigma_\infty\equiv0\) and \(\rho_\infty^{\mathrm{ac}}(b)=0\).
\end{proof}

We now prove Theorem~\ref{thm-stahl-limit}.

\begin{proof}
\emph{Atoms and nodes.} \(A_p^{-1/2}=p^{-1/2}(I+A_0/p)^{-1/2}\) gives
\(L_p=\widetilde B_p/p\) with
\(\widetilde B_p:=(I+A_0/p)^{-1/2}B(I+A_0/p)^{-1/2}\to B\), so by Kato
perturbation (\(B\) simple) \(p\,l_i^{(p)}\to b_i\). Moreover,
\(p\,N_i^{(p)}=1-c_i^{(p)}/p+O(p^{-2})\) with \(c_i^{(p)}\to a_i\), so
the rescaled atoms
\(p^{p}(N_i^{(p)})^{p}=e^{-c_i^{(p)}+O(1/p)}\to e^{-a_i}\):
\(\widetilde\rho_p^{\mathrm{at}}\to\sum_ie^{-a_i}\delta_{b_i}\) weakly.

\begin{enumerate}
\def\labelenumi{(\roman{enumi})}
\item
  Each \(\widetilde\rho_p\) is positive, and since
  \(p\,l_i^{(p)}\to b_i\) its support lies in the fixed compact
  \(J:=[b_1-1,b_n+1]\) for \(p\) large, with mass
  \(\to\operatorname{tr}e^{-A_0}\) by Lemma~\ref{lem-HY}. By weak-\(*\)
  compactness every sequence \(p_k\to\infty\) has a subsequence along
  which \(\widetilde\rho_{p_k}\to\nu\) weakly, \(\nu\ge0\). Since
  \((1+tb/p)^{-p}\to e^{-tb}\) uniformly on \(J\), Lemma~\ref{lem-HY}
  gives \(\int e^{-tb}\,d\nu=\operatorname{tr}\,e^{-(A_0+tB)}\) for all
  \(t\ge0\). A finite measure on a compact set is determined by its
  Laplace transform, so \(\nu=\rho_\infty\) for every subsequence.
\item
  is Lemma~\ref{lem-cts}.
\item
  By the first step and (i),
  \(\widetilde\rho_p^{\mathrm{ac}}=\widetilde\rho_p-\widetilde\rho_p^{\mathrm{at}}\to
  \nu^{\mathrm{c}}:=\rho_\infty-\sum_ie^{-a_i}\delta_{b_i}\) weakly, and
  \(\nu^{\mathrm c}\ge0\) as a weak limit of positive measures. \emph{No
  leak.} For \(h\ge0\) continuous, Fatou and (ii) give
  \(\int h\,d\nu^{\mathrm c}\ge\int h\,\rho_\infty^{\mathrm{ac}}\,db\),
  and on compact subsets of \((b_1,b_n)\setminus\{b_i\}\) the uniform
  bound of Lemma~\ref{lem-cts} upgrades this to equality. Hence the
  defect
  \(\varsigma:=\nu^{\mathrm c}-\rho_\infty^{\mathrm{ac}}\,db\ge0\) is
  carried by the nodes \(\{b_1,\dots,b_n\}\). It remains to see that
  \(\rho_\infty(\{b_i\})=e^{-a_i}\), so that \(\nu^{\mathrm c}\) is
  atomless and \(\varsigma=0\). Both sides of (\ref{eq-stahl-laplace})
  are entire in \(t\) (\(\rho_\infty\) is a finite measure on a compact
  set), so
  \(\hat\rho_\infty(s):=\int e^{-isb}\,d\rho_\infty(b)=\operatorname{tr}e^{-(A_0+isB)}\)
  for \(s\in\mathbb R\). For \(|s|\) large,
  \(A_0+isB=is(B-\tfrac{i}{s}A_0)\) and, \(B\) having simple spectrum,
  the eigenvalues of \(B-\tfrac isA_0\) are
  \(b_j-\tfrac isa_j+O(s^{-2})\) (\citeproc{ref-kato1995}{Kato 1995},
  Thm. II.5.11). Hence
  \(\hat\rho_\infty(s)=\sum_je^{-isb_j}e^{-a_j}+O(1/|s|)\). By Wiener's
  theorem,
  \(\rho_\infty(\{b_i\})=\lim_{T\to\infty}\tfrac1{2T}\int_{-T}^{T}e^{isb_i}\hat\rho_\infty(s)\,ds\),
  and the terms \(j\neq i\) and the remainder have Cesàro mean zero,
  leaving \(e^{-a_i}\). (This is the residue computation of Stahl, as
  presented by Eremenko (\citeproc{ref-eremenko2015}{Eremenko 2015},
  Prop. 1), in Fourier form.) Hence
  \(\rho_\infty=\sum_ie^{-a_i}\delta_{b_i}+\rho_\infty^{\mathrm{ac}}\,db\).
\end{enumerate}

\end{proof}

\subsection{\texorpdfstring{The \(n=2\) case: Bessis--Moussa--Villani's
Bessel
density}{The n=2 case: Bessis--Moussa--Villani's Bessel density}}\label{sec-n2-stahl}

For \(n=2\) the density (\ref{eq-phi-def}) is a Bessel function.

\begin{proposition}[Stahl density,
\(n=2\)]\protect\hypertarget{prp-n2-stahl}{}\label{prp-n2-stahl}

Let \(n=2\) and, in the \(B\)-eigenbasis,
\(B=\operatorname{diag}(b_1,b_2)\) with \(b_1<b_2\) and
\(A_0=\bigl(\begin{smallmatrix}a_1&c\\\bar c&a_2\end{smallmatrix}\bigr)\).
The Stahl measure of \(\operatorname{tr}\,e^{-(A_0+tB)}\) is
\(e^{-a_1}\delta_{b_1}+e^{-a_2}\delta_{b_2}+\rho_\infty^{\mathrm{ac}}(b)\,db\)
with, on \((b_1,b_2)\),
\begin{equation}\phantomsection\label{eq-n2-density}{
\begin{aligned}
\rho_\infty^{\mathrm{ac}}(b)&=\frac{|c|}{\sqrt{(b-b_1)(b_2-b)}}\,e^{\tau_\star(b)}\,I_1\bigl(\tau_w(b)\bigr),\\
\tau_\star(b)&=-\frac{a_1(b_2-b)+a_2(b-b_1)}{b_2-b_1},\qquad
\tau_w(b)=\frac{2|c|\sqrt{(b-b_1)(b_2-b)}}{b_2-b_1},
\end{aligned}
}\end{equation} \(I_1\) the modified Bessel function. The density is
bounded, with
\(\rho_\infty^{\mathrm{ac}}(b_1^+)=|c|^2e^{-a_1}/(b_2-b_1)\), and is
\(O(|c|^2)\) as \([A_0,B]\to0\).

\end{proposition}

\begin{proof}
In the \(B\)-eigenbasis, with \(\alpha_i:=a_i+\tau\) and
\(\beta_i:=b_i-b\) (\(\beta_1<0<\beta_2\)),
\(K_\infty(\tau,b)=(A_0+\tau I)(B-bI)^{-1}\) has trace
\((\alpha_1\beta_2+\alpha_2\beta_1)/(\beta_1\beta_2)\) and determinant
\((\alpha_1\alpha_2-|c|^2)/(\beta_1\beta_2)\), so its eigenvalues are a
conjugate pair with
\(\Sigma_\infty=\sqrt{4|\beta_1\beta_2||c|^2-(\alpha_1\beta_2-\alpha_2\beta_1)^2}\,/|\beta_1\beta_2|\)
exactly when the radicand is positive. Since
\(\alpha_1\beta_2-\alpha_2\beta_1=(b_2-b_1)(\tau-\tau_\star(b))\) and
\(|\beta_1\beta_2|=(b-b_1)(b_2-b)\), this is the window
\(|\tau-\tau_\star|<\tau_w\), on which
\(\Sigma_\infty=(b_2-b_1)\sqrt{\tau_w^2-(\tau-\tau_\star)^2}/|\beta_1\beta_2|\).
Inserting into (\ref{eq-phi-def}) with \(\tau=\tau_\star+s\) and
\(\int_{-w}^{w}e^{s}\sqrt{w^2-s^2}\,ds=\pi wI_1(w)\) gives
(\ref{eq-n2-density}). The atoms are those of
Theorem~\ref{thm-stahl-limit}(iii). Boundedness at the edges is
\(I_1(\tau_w)\sim\tau_w/2\).
\end{proof}

\begin{figure}[tb]

\centering{

\includegraphics[width=0.62\linewidth,height=\textheight,keepaspectratio]{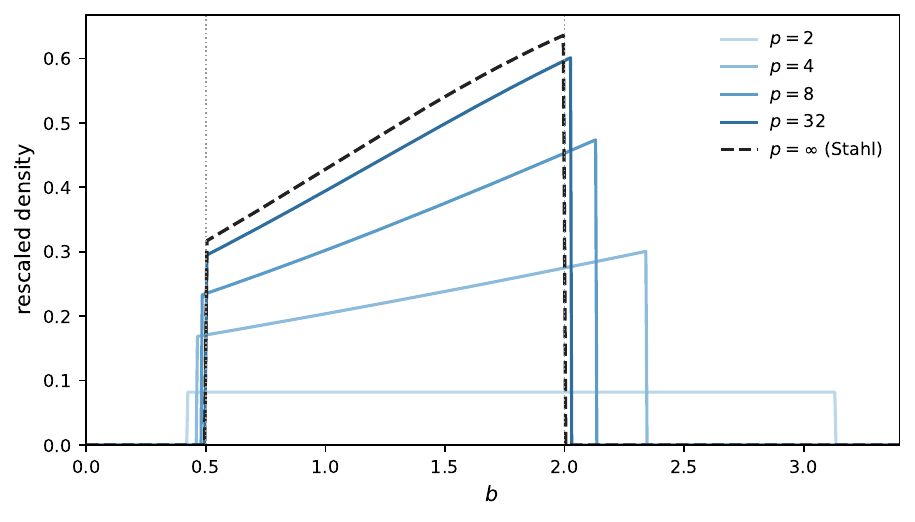}

}

\caption{\label{fig-stahl-confluence}The rescaled continuous densities
\(\widetilde\rho_p^{\mathrm{ac}}\) of Theorem~\ref{thm-stahl-limit} at
\(p=2,4,8,32\) for a \(2\times2\) pair
(\(B=\operatorname{diag}(0.5,2)\), and \(A_0\) with diagonal entries
\(0.3,-0.4\) and off-diagonal entry \(0.8\)), and their limit, the
Bessel density of Stahl's measure given in
Proposition~\ref{prp-n2-stahl} (dashed). The supports shrink to
\([b_1,b_2]\) (dotted lines). The atoms, not shown, converge to
\(e^{-a_i}\delta_{b_i}\). In this example the convergence is of order
\(1/p\).}

\end{figure}

This is the measure found by Bessis, Moussa and Villani
(\citeproc{ref-bessis1975}{Bessis, Moussa, and Villani 1975, sec. II}),
who describe it, for \(2\times2\) matrices, as a sum of delta functions
at the eigenvalues of \(B\) plus ``an entire and positive Bessel
function spread out between these eigenvalues''. Here it is recovered,
with its positivity, as the \(p\to\infty\) limit of the finite-\(p\)
family (Figure~\ref{fig-stahl-confluence}): along \(A_p=pI+A_0\) the
rescaled density \(p^{p-1}g_p(b/p)\), with \(g_p\) as in
(\ref{eq-n2p-density}), converges to the density (\ref{eq-n2-density}),
through
\({}_2F_1\bigl(\tfrac{1+p}{2},\tfrac{2+p}{2};2;\tfrac{w}{p^2}\bigr)\to
{}_0F_1(;2;\tfrac w4)=\tfrac{2}{\sqrt w}I_1(\sqrt w)\). For integer
\(p\) the finite-\(p\) densities are the B-splines of
Section~\ref{sec-bsplines}, and at \(p=2\) a box.

\section{\texorpdfstring{Order-one representation for
\(0<p<1\)}{Order-one representation for 0\textless p\textless1}}\label{sec-pless1}

For \(0<p<1\) the representation of Theorem~\ref{thm-rep} still holds,
but the prefactor \(p(p-1)\) is negative, so the order-\(p\) measure
\(\rho_p\) is \emph{signed}: positive atoms over a nonpositive
continuous part, positive as a whole if and only if \([A,B]=0\). Off the
commuting locus \(F_p\) is not a Stieltjes function of order \(p\). What
survives unconditionally is a positive representation of \emph{order
one}, a tracial Loewner representation of \(x\mapsto x^{-p}\), built
from the Stieltjes formula for \(x^{-p}\) and the \(p=1\) measures of
the shifted pairs \((A+uI,B)\). We construct it, give its density, and
relate it to \(\rho_p\).

On \((0,\infty)\), \(x\mapsto x^{-p}\) is a Stieltjes function (in
particular, its negative is a Pick function) with
\begin{equation}\phantomsection\label{eq-pick-formula}{
x^{-p}=\frac{\sin(p\pi)}{\pi}\int_0^\infty\frac{u^{-p}}{x+u}\,du,
}\end{equation} so by the Borel functional calculus and trace cyclicity
\begin{equation}\phantomsection\label{eq-F-via-pick}{
F_p(t)=\operatorname{tr}[(A+tB)^{-p}]
=\frac{\sin(p\pi)}{\pi}\int_0^\infty u^{-p}\operatorname{tr}[(A+uI+tB)^{-1}]\,du .
}\end{equation} The inner trace is the \(p=1\) case at \((A+uI,B)\). Put
\(L(u):=(A+uI)^{-1/2}B(A+uI)^{-1/2}\) (\(L(0)=L\)), with eigenpairs
\(l_i(u),e_i(u)\) and weights
\(N_i(u):=\langle e_i(u)|(A+uI)^{-1}|e_i(u)\rangle\). Then
\begin{equation}\phantomsection\label{eq-shifted-p1}{
\operatorname{tr}[(A+uI+tB)^{-1}]=\sum_{i=1}^n\frac{N_i(u)}{1+t\,l_i(u)} .
}\end{equation}

The argument is that of Hansen (\citeproc{ref-hansen2006trace}{Hansen
2006}, Main Theorem), who writes a nonnegative operator monotone
decreasing function \(f\) as a positive mixture of the functions
\(x\mapsto(c+x)^{-1}\) and concludes that
\(t\mapsto\operatorname{tr}f(A+tB)\) is the Laplace transform of a
positive measure. Lieb and Seiringer use the same device
(\citeproc{ref-lieb-seiringer2012}{Lieb and Seiringer 2012}, proof of
Thm. 2(b)). Here it is kept at the level of Stieltjes transforms.

\begin{proposition}[Order-\(1\) representation,
\(0<p<1\)]\protect\hypertarget{prp-stieltjes1}{}\label{prp-stieltjes1}

For \(A\succ0\), \(B\succeq0\), and \(0<p<1\), let \(\nu_u\) denote the
\(p=1\) spectral measure of the shifted pair \((A+uI,B)\) in
(\ref{eq-shifted-p1}). There is a unique positive finite Borel measure
\(\rho^{(1)}_{A,B,p}\) on \([0,l_n]\), characterized for every Borel set
\(E\subset[0,l_n]\) by
\begin{equation}\phantomsection\label{eq-rho1-mixture}{
\rho^{(1)}_{A,B,p}(E)
=\frac{\sin(\pi p)}{\pi}\int_0^\infty u^{-p}\nu_u(E)\,du,
}\end{equation} such that
\begin{equation}\phantomsection\label{eq-stieltjes-1-pless}{
F_p(t)=\int_0^{l_n}\frac{1}{1+t\ell}\,d\rho^{(1)}_{A,B,p}(\ell),\qquad t\ge0.
}\end{equation} Its total mass is \(\operatorname{tr}(A^{-p})\).

\end{proposition}

\begin{proof}
For continuous \(\varphi\) on \([0,l_n]\) the map
\(u\mapsto\int\varphi\,d\nu_u=\operatorname{tr}[\varphi(L(u))(A+uI)^{-1}]\)
is continuous, so \(u\mapsto\nu_u(E)\) is Borel for every Borel set
\(E\), and by Tonelli's theorem (\ref{eq-rho1-mixture}) defines a
positive measure. It is carried by \([0,l_n]\) because
\(\|L(u)\|=\|B^{1/2}(A+uI)^{-1}B^{1/2}\|\le\|B^{1/2}A^{-1}B^{1/2}\|=l_n\),
and its mass is \[
\frac{\sin(\pi p)}{\pi}\int_0^\infty u^{-p}\operatorname{tr}(A+uI)^{-1}\,du
=\operatorname{tr}(A^{-p}),
\] by (\ref{eq-pick-formula}) applied to \(A\). Substituting
(\ref{eq-shifted-p1}) into (\ref{eq-F-via-pick}) gives
(\ref{eq-stieltjes-1-pless}), and uniqueness is that of the Stieltjes
transform of a finite measure on a compact interval.
\end{proof}

The order-one measure is determined by the signed measure \(\rho_p\) of
Theorem~\ref{thm-rep}.

\begin{proposition}[Relation to
\(\rho_p\)]\protect\hypertarget{prp-bridge}{}\label{prp-bridge}

Let \(A\succ0\), \(B\succeq0\) and \(0<p<1\), and let \(\rho_p\) be the
signed measure of Theorem~\ref{thm-rep}. Then
\begin{equation}\phantomsection\label{eq-rho1-decomposition}{
\rho^{(1)}_{A,B,p}=w_0\,\delta_0+r_p(m)\,dm,
}\end{equation} where
\(w_0=\rho_p(\{0\})=\operatorname{tr}[(P_0A^{-1}P_0)^p]\), with \(P_0\)
the spectral projection of \(L\) onto \(\ker L\) (so \(w_0=0\) if
\(B\succ0\)), and, for almost every \(m\in(0,l_n)\),
\begin{equation}\phantomsection\label{eq-bridge}{
r_p(m)=\frac{\sin p\pi}{\pi}\,m^{p-1}\int_{(m,\,l_n]}(\ell-m)^{-p}\,d\rho_p(\ell).
}\end{equation} In particular \(\rho^{(1)}_{A,B,p}\) has no atoms on
\((0,l_n]\) and no singular continuous part.

\end{proposition}

\begin{proof}
For \(\ell>0\), the formula (\ref{eq-pick-formula}) at \(x=1+t\ell\),
with the substitution \(u=(\ell-m)/m\), reads
\((1+t\ell)^{-p}=\tfrac{\sin p\pi}{\pi}\int_0^\ell m^{p-1}(\ell-m)^{-p}(1+tm)^{-1}\,dm\).
Insert this into (\ref{eq-main-rep}) and exchange the integrals, which
is permitted because
\(\int\!\int_0^\ell m^{p-1}(\ell-m)^{-p}\,dm\,d|\rho_p|(\ell)=\tfrac{\pi}{\sin p\pi}\,|\rho_p|((0,l_n])\)
is finite. This exhibits \(F_p\) as the Stieltjes transform of the
finite signed measure \(\rho_p(\{0\})\,\delta_0+\tilde r(m)\,dm\), where
\(\tilde r\) is the right side of (\ref{eq-bridge}). A finite signed
measure on \([0,l_n]\) is determined by its Stieltjes transform, so this
measure is \(\rho^{(1)}_{A,B,p}\). The value of \(\rho_p(\{0\})\) is
given by (\ref{eq-atomic-Omega}).
\end{proof}

Formula (\ref{eq-bridge}) is the connection formula of Karp and
Prilepkina (\citeproc{ref-karp2012generalized}{Karp and Prilepkina
2012}, Thm. 4, (14)) between the representing measures of orders \(p\)
and \(1\), stated there for positive measures. In their terminology it
is a right-sided Kober--Erdélyi fractional integral of order \(1-p\).
Here it is applied to the signed measure \(\rho_p\), and the fractional
integration restores the positivity that \(\rho_p\) lacks off the
commuting locus.

The density also has a closed form. For \(\ell>0\) put \[
H_\ell:=\ell^{-1}B-A.
\]

\begin{theorem}[The order-\(1\)
density]\protect\hypertarget{thm-rho1-density}{}\label{thm-rho1-density}

Let \(A\succ0\), \(B\succeq0\), and \(0<p<1\). The density \(r_p\) of
\(\rho^{(1)}_{A,B,p}\) on \((0,l_n)\) is given, at every \(\ell\) that
is not an eigenvalue of \(L\), by
\begin{equation}\phantomsection\label{eq-rho-1-density}{
\boxed{\quad
r_p(\ell)
=\frac{\sin(\pi p)}{\pi\ell}\,
\operatorname{tr}\!\left[
\mathbf 1_{(0,\infty)}(H_\ell)H_\ell^{-p}
\right],\qquad H_\ell=\ell^{-1}B-A.\quad}
}\end{equation} In the open interval \((0,l_n)\) the only singularities
of \(r_p\) occur at the positive eigenvalues \(\ell_0\) of \(L\). If
\(E=\ker H_{\ell_0}\) and \(\beta_1,\ldots,\beta_k>0\) are the
eigenvalues of \(P_EBP_E\vert_E\), then
\begin{equation}\phantomsection\label{eq-rho1-threshold-asymptotic}{
r_p(\ell)
=\frac{\sin(\pi p)}{\pi\ell_0}
\sum_{j=1}^k\left(\frac{\beta_j}{\ell_0^2}\right)^{-p}
(\ell_0-\ell)^{-p}+O(1)
}\end{equation} as \(\ell\uparrow\ell_0\), and \(r_p\) is bounded as
\(\ell\downarrow\ell_0\). Thus every threshold singularity is one-sided
and locally integrable. At the left end,
\(r_p(\ell)\sim\tfrac{\sin\pi p}{\pi}\operatorname{tr}[B_+^{-p}]\,\ell^{p-1}\)
as \(\ell\downarrow0\), where \(B_+^{-p}\) is the power of \(B\) on its
range. This too is integrable.

\end{theorem}

\begin{proof}
\textbf{1. The density.} The family \(L(u)\), \(u\ge0\), is Hermitian
and real-analytic, so its eigenvalues \(l_i(u)\) and an orthonormal
basis of eigenvectors \(e_i(u)\), \(i=1,\dots,n\), can be chosen
analytic in \(u\), through crossings (Rellich; see
(\citeproc{ref-kato1995}{Kato 1995}, Thm. II.6.1 and §II.6.2)). Since
\(\tfrac{d}{du}(A+uI)^{-1/2}=-\tfrac12(A+uI)^{-3/2}\), \[
L'(u)=-\tfrac12\bigl[(A+uI)^{-1}L(u)+L(u)(A+uI)^{-1}\bigr],
\] and the Hellmann--Feynman formula gives
\begin{equation}\phantomsection\label{eq-hellmann-feynman}{
l_i'(u)=\langle e_i(u)|L'(u)|e_i(u)\rangle=-\,l_i(u)\,N_i(u).
}\end{equation} Hence \(l_i(u)=l_i(0)\exp\bigl(-\int_0^uN_i\bigr)\). A
branch with \(l_i(0)=0\) vanishes identically. A branch with
\(l_i(0)>0\) is positive and strictly decreasing, and it tends to \(0\)
as \(u\to\infty\) because \(N_i(u)\ge(\|A\|+u)^{-1}\). It is thus a
diffeomorphism of \((0,\infty)\) onto \((0,l_i(0))\), with inverse
\(u_i(\ell)\), and \(|d\ell|=\ell\,N_i(u)\,du\) along it. For a Borel
set \(E\subset(0,\infty)\), the mixture (\ref{eq-rho1-mixture}) and this
change of variables give \[
\rho^{(1)}_{A,B,p}(E)
=\frac{\sin(\pi p)}{\pi}\sum_{i:\,l_i(0)>0}\int_0^\infty u^{-p}N_i(u)\,\mathbf 1_E(l_i(u))\,du
=\frac{\sin(\pi p)}{\pi}\int_E\frac1\ell\sum_{i:\,l_i(0)>\ell}u_i(\ell)^{-p}\,d\ell .
\] Finally, for \(\ell>0\) and \(u\ge0\) the equation \(Bv=\ell(A+uI)v\)
says both that \(\ell\) is an eigenvalue of \(L(u)\), with eigenvector
\((A+uI)^{1/2}v\), and that \(u\) is an eigenvalue of \(H_\ell\), with
eigenvector \(v\). So the numbers \(u_i(\ell)\) with \(l_i(0)>\ell\) are
the positive eigenvalues of \(H_\ell\), counted with multiplicity, and
the last integrand is (\ref{eq-rho-1-density}), at every \(\ell\) that
is not an eigenvalue of \(L\).

\textbf{2. Thresholds.} The exceptional points of
(\ref{eq-rho-1-density}) are those \(\ell_0>0\) for which \(H_{\ell_0}\)
is singular, that is, \(Bv=\ell_0Av\) for some \(v\ne0\): they are the
positive eigenvalues of \(L\). Let \(E=\ker H_{\ell_0}\). Since
\(A\succ0\) and \(Bv=\ell_0Av\) on \(E\), the compression
\(P_EBP_E\vert_E\) is positive definite. Let its eigenvalues be
\(\beta_1,\ldots,\beta_k>0\). Moreover, \[
H_\ell-H_{\ell_0}
=(\ell^{-1}-\ell_0^{-1})B
=\frac{\ell_0-\ell}{\ell_0^2}B+O((\ell-\ell_0)^2).
\] Finite-dimensional degenerate perturbation theory
(\citeproc{ref-kato1995}{Kato 1995}, Thm. II.5.11 and Thm. II.6.1)
therefore gives the \(k\) eigenvalue branches emerging from zero as \[
u_j(\ell)
=\frac{\beta_j}{\ell_0^2}(\ell_0-\ell)+O((\ell-\ell_0)^2).
\] They are positive to the left of \(\ell_0\) and negative to the
right. Substitution into (\ref{eq-rho-1-density}) yields
(\ref{eq-rho1-threshold-asymptotic}). Since \(0<p<1\),
\((\ell_0-\ell)^{-p}\) is locally integrable. As \(\ell\downarrow0\) the
positive eigenvalues of \(H_\ell=\ell^{-1}(B-\ell A)\) are
\(b_j/\ell+O(1)\), with \(b_j\) the positive eigenvalues of \(B\), which
gives the behaviour at the left end.
\end{proof}

The atoms of \(\rho_p\) produce, through (\ref{eq-bridge}), the
threshold singularities (\ref{eq-rho1-threshold-asymptotic}): the atom
\(\Omega_s\) at \(\lambda_s=\ell_0\) contributes
\(\tfrac{\sin p\pi}{\pi}\,\ell_0^{p-1}\Omega_s\,(\ell_0-m)^{-p}(1+o(1))\)
as \(m\uparrow\ell_0\), which is the leading term there. Indeed
\(E=\ker(B-\ell_0A)=A^{-1/2}V_s\) and \(P_EBP_E|_E=\ell_0P_EAP_E|_E\),
whose eigenvalues are \(\beta_j=\ell_0/\nu_{s,j}\) by the computation in
the proof of Lemma~\ref{lem-Apole}(i). Hence
\(\ell_0^{-1}\sum_j(\beta_j/\ell_0^2)^{-p}=\ell_0^{p-1}\sum_j\nu_{s,j}^{\,p}=\ell_0^{p-1}\Omega_s\).

\appendix

\section{The connection with quantum distributionally robust
optimization}\label{sec-appendix-dro}

This note grew out of the companion paper
(\citeproc{ref-reybellet2026dro}{Rey-Bellet 2026}), which uses several
of its statements. We gather them here, with the dictionary between the
two notations.

\textbf{The setting.} The companion paper studies the largest
expectation \(\operatorname{tr}(L\rho)\) of an observable \(L=L^*\) over
the states \(\rho\) in a ball around a faithful state \(\sigma\), the
ball being defined by the sandwiched \(\alpha\)-divergence
\(H_\alpha(\rho\|\sigma)=\frac{1}{\alpha-1}\bigl(\operatorname{tr}[(\sigma^{\frac{1-\alpha}{2\alpha}}\rho\,\sigma^{\frac{1-\alpha}{2\alpha}})^{\alpha}]-1\bigr)\)
(the sandwiched Rényi divergence in its power normalization), of order
\(\alpha\ge\tfrac12\), \(\alpha\neq1\). With
\(\beta:=(\alpha-1)/(2\alpha)\), the dual of this problem is expressed,
for \(\tfrac12\le\alpha<1\), through the trace function \[
c\mapsto\operatorname{tr}\bigl[\bigl(\sigma^\beta(cI-L)\sigma^\beta\bigr)^{-p}\bigr],\qquad p=\frac{\alpha}{1-\alpha}.
\] (For \(\alpha>1\) the corresponding function is
\(\operatorname{tr}[(\sigma^\beta(L-cI)\sigma^\beta)_+^{\alpha/(\alpha-1)}]\).
It is covered directly by Heinävaara's identity and is treated in the
companion paper.) Writing it as a transform of a positive measure
\(\omega_\alpha\) on the spectral interval of \(L\) turns the quantum
problem into a classical one, and the concentration inequalities and
Padé-type certificates of (\citeproc{ref-reybellet2026dro}{Rey-Bellet
2026}) rest on this step. Such certificates were, not coincidentally,
the original motivation of Bessis, Moussa and Villani
(\citeproc{ref-bessis1975}{Bessis, Moussa, and Villani 1975}). In the
notation of Corollary~\ref{cor-pencil} below, the trace function is
\(\operatorname{tr}[(cP-Q)^{-p}]\) with \(P=\sigma^{2\beta}\) and
\(Q=\sigma^\beta L\sigma^\beta\), so that \(P^{-1/2}QP^{-1/2}\) is the
observable \(L\) itself.

\textbf{What is used.} The companion paper imports Theorem~\ref{thm-rep}
itself, with the density (\ref{eq-cts}). For the support of
\(\mu_{A,B}\) it uses the cone statement of Lemma~\ref{lem-cone}(ii)
and, for the integrability of the radial weight, the bound
(\ref{eq-support-bound}) obtained in the proof of Theorem~\ref{thm-rep}.
The inclusion of \(\operatorname{supp}\mu_c\) in the joint numerical
range is Heinävaara's (\citeproc{ref-heinavaara-thesis}{Heinävaara
2024b}, Prop. 4) and is not needed. For orders \(\tfrac12\le\alpha<1\):
the pencil form of the representation, with its atoms at every
multiplicity (Corollary~\ref{cor-pencil}), specialized below. For
integer \(p\): the B-spline form (Theorem~\ref{thm-bspline}), again
specialized below (\(p=2\) is \(\alpha=\tfrac23\)). For \(\alpha=1\):
the atoms and the density of Stahl's measure
(Theorem~\ref{thm-stahl-limit}). For a qubit: the closed forms of
Section~\ref{sec-n2}.

\subsection{Orders below one: the pencil form}\label{sec-pencil-form}

The representation extends to a Hermitian pencil, as a generalized
Stieltjes transform in the spectral variable of \(P^{-1/2}QP^{-1/2}\): a
shift and a Möbius change of variable deliver it, at every multiplicity.
This is the form used in the companion paper.

\begin{corollary}[Pencil
form]\protect\hypertarget{cor-pencil}{}\label{cor-pencil}

Let \(P\succ0\) and \(Q\) be Hermitian on \(\mathbb C^n\),
\(\Lambda:=P^{-1/2}QP^{-1/2}\) with distinct eigenvalues
\(m_1<\dots<m_r\) and spectral projections \(\Pi_1,\dots,\Pi_r\), and
\(p>0\). There is a unique finite signed measure \(\omega\) on
\([m_1,m_r]\) with
\begin{equation}\phantomsection\label{eq-pencil-form}{
\operatorname{tr}\bigl[(cP-Q)^{-p}\bigr]=\int(c-m)^{-p}\,d\omega(m),\qquad c>m_r ,
}\end{equation} namely
\(\omega=\mathfrak m_*\bigl(\ell'^{\,-p}\rho_p\bigr)\), where \(\rho_p\)
is the measure of Theorem~\ref{thm-rep} for the pair \((c_0P-Q,\,P)\),
\(c_0>m_r\) arbitrary, and \(\mathfrak m(\ell'):=c_0-1/\ell'\). Its mass
is \(\operatorname{tr}P^{-p}\), its atoms are
\begin{equation}\phantomsection\label{eq-pencil-atoms}{
\omega(\{m_s\})=\operatorname{tr}\bigl[(\Pi_sP^{-1}\Pi_s)^{p}\bigr],
}\end{equation} its continuous part vanishes iff \(p=1\) or \([P,Q]=0\),
and \(\omega\ge0\) for \(p\ge1\).

\end{corollary}

\begin{proof}
With \(A:=c_0P-Q\succ0\), \(B:=P\succ0\) and \(t:=c-c_0\),
\(A+tB=cP-Q\), and \(A^{-1}B\) is similar to \((c_0I-\Lambda)^{-1}\), so
the spectral interval of the pair is \([\ell'_1,\ell'_r]\),
\(\ell'_s=(c_0-m_s)^{-1}\), which \(\mathfrak m\) maps increasingly onto
\([m_1,m_r]\). Since \(1+t\ell'=(c-m)/(c_0-m)\) for
\(m=\mathfrak m(\ell')\),
i.e.~\((1+t\ell')^{-p}=(c-m)^{-p}\ell'^{\,-p}\), (\ref{eq-main-rep})
pushed forward through \(\mathfrak m\) is (\ref{eq-pencil-form}) for
\(c\ge c_0\), and for every \(c>m_r\) by the analytic continuation of
(\ref{eq-main-rep}) to \(t>-1/\ell'_r=m_r-c_0\) established in the proof
of Theorem~\ref{thm-rep}. The mass is the limit
\(c^p\operatorname{tr}(cP-Q)^{-p}\to
\operatorname{tr}P^{-p}\), and uniqueness follows from the moments, as
in Theorem~\ref{thm-rep}. For the atoms, the \(\ell'_s\)-eigenspace of
\(A^{-1/2}BA^{-1/2}\) is \(A^{1/2}P^{-1/2}V_s\) with \(V_s\) the
\(m_s\)-eigenspace of \(\Lambda\): if \(u_1,\dots,u_{k_s}\) is an
orthonormal basis of \(V_s\) then
\(v_a:=A^{1/2}P^{-1/2}u_a/\sqrt{c_0-m_s}\) is an orthonormal basis of
that eigenspace, since
\(\langle P^{-1/2}u_a|A|P^{-1/2}u_b\rangle=\langle u_a|c_0I-\Lambda|u_b\rangle=(c_0-m_s)\delta_{ab}\),
and the compression of \(A^{-1}\) to it has matrix
\(\langle v_a|A^{-1}|v_b\rangle=\langle u_a|P^{-1}|u_b\rangle/(c_0-m_s)\).
Hence the \(\rho_p\)-atom
\(\operatorname{tr}[(P_s^{\rm pair}A^{-1}P_s^{\rm pair})^p]=(c_0-m_s)^{-p}\operatorname{tr}[(\Pi_sP^{-1}\Pi_s)^p]\),
and the factor \(\ell'^{\,-p}_s=(c_0-m_s)^{p}\) carried by the
pushforward gives (\ref{eq-pencil-atoms}). The continuous part is the
image of \(\ell'^{\,-p}\rho_p^{\mathrm{ac}}\) and vanishes iff \(p=1\)
or \([A,B]=[P,Q]=0\). Positivity for \(p\ge1\) is that of \(\rho_p\).
\end{proof}

In Corollary~\ref{cor-pencil} take, as in the companion paper
(\citeproc{ref-reybellet2026dro}{Rey-Bellet 2026, sec. 4}),
\(P=\sigma^{2\beta}\) and \(Q=\sigma^\beta L\sigma^\beta\) with
\(\beta=(\alpha-1)/(2\alpha)<0\) for a Rényi order
\(\alpha\in[\tfrac12,1)\), so that \(\Lambda=L\) and \(\Pi_s=P_s\),
\(p=\alpha/(1-\alpha)\) and
\(\operatorname{tr}P^{-p}=\operatorname{tr}\sigma=1\): \(\omega\) is the
probability measure \(\omega_\alpha\) that classicalizes the quantum DRO
problem, with atoms
\(\operatorname{tr}[(P_s\sigma^{-2\beta}P_s)^{\alpha/(1-\alpha)}]\) at
the eigenvalues of \(L\), whatever their multiplicities. (The paper's
\(N_i=\langle i|\sigma^{-2\beta}|i\rangle\) is the present \(N_i\) of
the pair \((c_0P-Q,P)\) multiplied by \(c_0-m_i\), as the pushforward
shows.)

For integer \(p\) the expansion of Section~\ref{sec-bsplines} applies to
the pencil form (\ref{eq-pencil-form}): for \(P\succ0\), \(Q\) Hermitian
and integer \(p\ge2\),
\(\operatorname{tr}[(cP-Q)^{-p}]=\operatorname{tr}[((cI-\Lambda)^{-1}P^{-1})^p]\)
in the eigenbasis of \(\Lambda=P^{-1/2}QP^{-1/2}\) is
(\ref{eq-Bsplines}) with \((c-m_i)^{-1}\) in place of \((1+tl_i)^{-1}\)
and \(N_{ij}=\langle e_i|P^{-1}|e_j\rangle\), and the Peano identity
holds for the kernel \((c-m)^{-p}\) as it does for \((1+t\ell)^{-p}\).
So the measure \(\omega\) of Corollary~\ref{cor-pencil} is
\(\sum N_{i_1i_2}\cdots N_{i_pi_1}M(\cdot\mid m_{i_1},\dots,m_{i_p})\)
with knots at \(\operatorname{spec}\Lambda\). At \(p=2\), for every
\(n\) and when \(\Lambda\) has simple spectrum,
\begin{equation}\phantomsection\label{eq-pencil-p2}{
\omega=\sum_iN_{ii}^{\,2}\,\delta_{m_i}+\sum_{i<j}\frac{2|N_{ij}|^2}{m_j-m_i}\,\mathbf 1_{[m_i,m_j]}\,dm ,
}\end{equation} which with \(P=\sigma^{2\beta}\) is the measure
\(\omega_{2/3}\) of the companion paper.

\section*{Acknowledgments}\label{acknowledgments}
\addcontentsline{toc}{section}{Acknowledgments}

This note grew out of work on quantum distributionally robust
optimization (\citeproc{ref-reybellet2026dro}{Rey-Bellet 2026}),
supported by the National Science Foundation under grants DMS-2307115
and DMS-2606221. I thank the Département de Mathématiques of the École
Polytechnique Fédérale de Lausanne for its hospitality while this work
was completed. Throughout this project I used Claude (Anthropic) as a
research assistant: to help develop some of the arguments, to search the
literature and cross-check proofs and references, to draft and
restructure the text, and to write the code that produces the figures. I
verified all mathematical statements and take full responsibility for
the content.

\section*{References}\label{references}
\addcontentsline{toc}{section}{References}

\phantomsection\label{refs}

\bigskip
\noindent\textsc{Luc Rey-Bellet}\\
Department of Mathematics and Statistics, University of Massachusetts Amherst,\\
710 N.~Pleasant St., Amherst, MA 01003, USA\\
\textit{E-mail address:} \texttt{luc@cns.umass.edu}


\begin{CSLReferences}{1}{0}
\bibitem[\citeproctext]{ref-berta2017}
Berta, M., O. Fawzi, and M. Tomamichel. 2017. {``On Variational
Expressions for Quantum Relative Entropies.''} \emph{Lett. Math. Phys.}
107: 2239--65. \url{https://doi.org/10.1007/s11005-017-0990-7}.

\bibitem[\citeproctext]{ref-bessis1975}
Bessis, D., P. Moussa, and M. Villani. 1975. {``Monotonic Converging
Variational Approximations to the Functional Integrals in Quantum
Statistical Mechanics.''} \emph{J. Math. Phys.} 16: 2318--25.
\url{https://doi.org/10.1063/1.522463}.

\bibitem[\citeproctext]{ref-bhatia1997}
Bhatia, R. 1997. \emph{Matrix Analysis}. Graduate Texts in Mathematics
169. Springer. \url{https://doi.org/10.1007/978-1-4612-0653-8}.

\bibitem[\citeproctext]{ref-bullen1971}
Bullen, P. S. 1971. {``A Criterion for \(n\)-Convexity.''} \emph{Pacific
J. Math.} 36 (1): 81--98. \url{https://doi.org/10.2140/pjm.1971.36.81}.

\bibitem[\citeproctext]{ref-carlen2010}
Carlen, E. A. 2010. {``Trace Inequalities and Quantum Entropy: An
Introductory Course.''} In \emph{Entropy and the Quantum}, edited by R.
Sims and D. Ueltschi, 529:73--140. Contemp. Math. Amer. Math. Soc.
\url{https://doi.org/10.1090/conm/529/10428}.

\bibitem[\citeproctext]{ref-clivaz2016}
Clivaz, Fabien. 2016. {``Stahl's Theorem (Aka {BMV} Conjecture):
{Insights} and Intuition on Its Proof.''} In \emph{Spectral Theory and
Mathematical Physics}, edited by Marius Mantoiu, Georgi Raikov, and
Rafael Tiedra de Aldecoa, 107--17. Cham: Springer International
Publishing. \url{https://doi.org/10.1007/978-3-319-29992-1_6}.

\bibitem[\citeproctext]{ref-curry-schoenberg1966}
Curry, H. B., and I. J. Schoenberg. 1966. {``On {P}{ó}lya Frequency
Functions {IV}: The Fundamental Spline Functions and Their Limits.''}
\emph{J. Analyse Math.} 17: 71--107.
\url{https://doi.org/10.1007/BF02788653}.

\bibitem[\citeproctext]{ref-deBoor2001}
de Boor, Carl. 2001. \emph{A Practical Guide to Splines}. Revised. Vol.
27. Applied Mathematical Sciences. New York: Springer.

\bibitem[\citeproctext]{ref-deboor2005divided}
---------. 2005. {``Divided Differences.''} \emph{Surv. Approx. Theory}
1: 46--69. \url{https://arxiv.org/abs/math/0502036}.

\bibitem[\citeproctext]{ref-dykema-skripka2009}
Dykema, Ken, and Anna Skripka. 2009. {``Higher Order Spectral Shift.''}
\emph{J. Funct. Anal.} 257: 1092--132.
\url{https://doi.org/10.1016/j.jfa.2009.02.019}.

\bibitem[\citeproctext]{ref-eremenko2015}
Eremenko, A. 2015. {``Herbert {S}tahl's Proof of the {BMV}
Conjecture.''} \emph{Sb. Math.} 206: 87--92.
\url{https://doi.org/10.1070/SM2015v206n01ABEH004447}.

\bibitem[\citeproctext]{ref-frank-lieb2013}
Frank, R. L., and E. H. Lieb. 2013. {``Monotonicity of a Relative
{R}ényi Entropy.''} \emph{J. Math. Phys.} 54: 122201.
\url{https://doi.org/10.1063/1.4838835}.

\bibitem[\citeproctext]{ref-gradshteyn-ryzhik2015}
Gradshteyn, I. S., and I. M. Ryzhik. 2015. \emph{Table of Integrals,
Series, and Products}. 8th ed. Academic Press.
\url{https://doi.org/10.1016/C2010-0-64839-5}.

\bibitem[\citeproctext]{ref-hansen2006trace}
Hansen, Frank. 2006. {``Trace Functions as {L}aplace Transforms.''}
\emph{J. Math. Phys.} 47: 043504.
\url{https://doi.org/10.1063/1.2186925}.

\bibitem[\citeproctext]{ref-heinavaara2024planes}
Heinävaara, Otte. 2024a. {``Planes in {S}chatten-3.''} \emph{J. Funct.
Anal.} 287 (2): 110469. \url{https://doi.org/10.1016/j.jfa.2024.110469}.

\bibitem[\citeproctext]{ref-heinavaara-thesis}
---------. 2024b. {``Tracial Joint Spectral Measures.''} PhD thesis,
Princeton, NJ: Princeton University.
\url{https://www.its.caltech.edu/~oeh/thesis_final.pdf}.

\bibitem[\citeproctext]{ref-heinavaara2024tracial}
---------. 2025. {``Tracial Joint Spectral Measures.''} \emph{Invent.
Math.} 239 (2): 505--26.
\url{https://doi.org/10.1007/s00222-024-01309-6}.

\bibitem[\citeproctext]{ref-hiai2010}
Hiai, Fumio. 2010. {``Matrix Analysis: Matrix Monotone Functions, Matrix
Means, and Majorization.''} \emph{Interdiscip. Inform. Sci.} 16 (2):
139--246. \url{https://doi.org/10.4036/iis.2010.139}.

\bibitem[\citeproctext]{ref-Karlin1968}
Karlin, Samuel. 1968. \emph{Total Positivity, Volume {I}}. Stanford, CA:
Stanford University Press.

\bibitem[\citeproctext]{ref-karp2012generalized}
Karp, Dmitrii, and Elena Prilepkina. 2012. {``Generalized {S}tieltjes
Functions and Their Exact Order.''} \emph{J. Class. Anal.} 1 (1):
53--74. \url{https://doi.org/10.7153/jca-01-07}.

\bibitem[\citeproctext]{ref-kato1995}
Kato, T. 1995. \emph{Perturbation Theory for Linear Operators}. 2nd ed.
Springer. \url{https://doi.org/10.1007/978-3-642-66282-9}.

\bibitem[\citeproctext]{ref-lieb-seiringer2004}
Lieb, E. H., and R. Seiringer. 2004. {``Equivalent Forms of the
{B}essis-{M}oussa-{V}illani Conjecture.''} \emph{J. Stat. Phys.} 115:
185--90. \url{https://doi.org/10.1023/B:JOSS.0000019811.15510.27}.

\bibitem[\citeproctext]{ref-lieb-seiringer2012}
---------. 2012. {``Further Implications of the
{B}essis-{M}oussa-{V}illani Conjecture.''} \emph{J. Stat. Phys.} 149:
86--91. \url{https://doi.org/10.1007/s10955-012-0585-8}.

\bibitem[\citeproctext]{ref-pradhan-skripka2025}
Pradhan, Chandan, and Anna Skripka. 2025. {``Positivity of Spectral
Shift Functions and Infinite-Dimensional {BMV} Conjecture.''}
\url{https://doi.org/10.48550/arXiv.2512.05587}.

\bibitem[\citeproctext]{ref-reybellet2026dro}
Rey-Bellet, Luc. 2026. {``Quantum Robustness via Sandwiched {Rényi}
Divergences.''} Preprint, in preparation.

\bibitem[\citeproctext]{ref-skripka2011}
Skripka, Anna. 2011. {``Multiple Operator Integrals and Spectral
Shift.''} \emph{Illinois J. Math.} 55 (1): 305--24.
\url{https://doi.org/10.1215/ijm/1355927038}.

\bibitem[\citeproctext]{ref-sokal2010real}
Sokal, Alan D. 2010. {``Real-Variables Characterization of Generalized
{S}tieltjes Functions.''} \emph{Expo. Math.} 28 (2): 179--85.
\url{https://doi.org/10.1016/j.exmath.2009.06.004}.

\bibitem[\citeproctext]{ref-stahl2013}
Stahl, H. R. 2013. {``Proof of the {BMV} Conjecture.''} \emph{Acta
Math.} 211: 255--90. \url{https://doi.org/10.1007/s11511-013-0104-z}.

\bibitem[\citeproctext]{ref-takahashi-fujiwara2017}
Takahashi, Kaito, and Akio Fujiwara. 2017. {``Information Geometry of
Sandwiched {R{é}nyi} \(\alpha\)-Divergence.''} \emph{J. Phys. A} 50
(16): 165301. \url{https://doi.org/10.1088/1751-8121/aa6326}.

\bibitem[\citeproctext]{ref-widder1938stieltjes}
Widder, D. V. 1938. {``The {S}tieltjes Transform.''} \emph{Trans. Amer.
Math. Soc.} 43 (1): 7--60. \url{https://doi.org/10.2307/1989901}.

\end{CSLReferences}
\end{document}